\documentclass[12pt,a4paper]{article}

\usepackage[utf8]{inputenc}
\usepackage[T1]{fontenc}
\usepackage{lmodern}
\usepackage{microtype}
\usepackage{enumerate}

\usepackage[english]{babel}
\usepackage{booktabs}

\usepackage{amssymb}
\usepackage{mathtools}
\mathtoolsset{
  showonlyrefs,
  mathic
}

\usepackage[dvipsnames]{xcolor}

\usepackage{hyperref}
\hypersetup{%
  colorlinks=true,
  breaklinks=true,
  linkcolor=red,
  anchorcolor=black,
  citecolor=brown,
  filecolor=blue,
  menucolor=red,
  urlcolor=red,
}

\usepackage[round,longnamesfirst]{natbib}

\usepackage{graphicx}
\graphicspath{{img/}}

\makeatletter
\def\calFactory#1{%
   \expandafter\def\csname c#1\endcsname{\mathcal{#1}}}
\def\frakFactory#1{%
   \expandafter\def\csname k#1\endcsname{\mathfrak{#1}}}
\def\bbFactory#1{%
   \expandafter\def\csname b#1\endcsname{\mathbb{#1}}}
\def\bfFactory#1{%
   \expandafter\def\csname f#1\endcsname{\mathbf{#1}}}
\newcounter{ctr}
\loop
  \stepcounter{ctr}
  \edef\X{\@Alph\c@ctr}%
  \edef\x{\@alph\c@ctr}%
  \expandafter\calFactory\X
  \expandafter\frakFactory\X
  \expandafter\frakFactory\x
  \expandafter\bbFactory\X
  \expandafter\bfFactory\X
\ifnum\thectr<26
\repeat
\makeatother

\newcommand{\mathset}[1]{\mathbb{#1}}

\newcommand{\N}{\mathset{N}}

\newcommand{\R}{\mathset{R}}

\newcommand{\Leb}{\mathset{L}}
\newcommand{\Lds}{{\fL^2_{\mathrm{sym}}}}

\newcommand{\mathproba}[1]{\mathbf{#1}}

\newcommand{\e}{\mathproba{E}}
\newcommand{\p}{\mathproba{P}}
\newcommand{\1}{\mathproba{I}}

\DeclareMathOperator*{\argmax}{arg\,max}

\let\var\Var

\newcommand{\floor}[1]{\lfloor#1\rfloor}
\newcommand{\kkk}{\mathbf{k}}
\newcommand{\D}{\mathrm{d}}
\newcommand{\Klh}[2]{K^{(#1)}_{#2}}
\newcommand{\smax}{{s^*}}
\newcommand{\smin}{{s_*}}

\newtheorem{theorem}{Theorem}

\newtheorem{lemma}{Lemma}

\newtheorem{definition}{Definition}

\newtheorem{remark}{Remark}

\begin{document}

\title{Adaptive regression with Brownian path covariate}
\author{%
Karine~Bertin%
\thanks{%
CIMFAV-INGEMAT, Universidad de Valpa{r}a\'{\i}so, General Cruz 222,
Valparaíso, Chile,
\url{karine.bertin@uv.cl}}
\and
Nicolas~Klutchnikoff%
\thanks{%
Univ Rennes, CNRS, IRMAR~--~UMR 6625, F-35000 Rennes, France, 
\url{nicolas.klutchnikoff@univ-rennes2.fr}}
}

\date{}
\maketitle

\begin{abstract}
\noindent
This paper deals with estimation with functional covariates. More precisely, we aim at estimating the regression function $m$ of a continuous outcome $Y$ against a standard Wiener coprocess $W$. Following~\citet{MR3412645} and~\citet{MR3716123} the Wiener-Itô decomposition of $m(W)$ is used to construct a family of estimators. The minimax rate of convergence over specific smoothness classes is obtained. A data-driven selection procedure is defined following the ideas developed by \citet{GL2011}. An oracle-type inequality is obtained which leads to adaptive results.


\end{abstract}

\section{Introduction}

The problem of regression estimation  is one of the most studied in statistics and different models have been considered depending on the nature of the data.
In an increasing number of applications, it seems natural to assume that the covariate takes values in a functional space. The book of \citet{MR2168993} provides an overview on the subject of functional data analysis. In this context, several authors studied linear functional regression models \citep[see for example][]{MR2163159,MR2291496,MR2488344}. Nonparametric functional regression models have also been investigated \citep[see][and references therein]{MR2229687}.
In this paper, we are interested in such a model where the covariate is a Wiener Process.
More precisely, let $\varepsilon$ be a real-valued random variable and $W=(W(t) : 0\leq t \leq 1)$ be a standard Brownian motion independent of $\varepsilon$. We define
\begin{equation}\label{eq:model}
  Y = m(W) + \varepsilon
\end{equation}
where $m:\cC \to \R$ is a mapping defined on the set $\cC$ of all continuous functions $w:[0,1]\to\R$ and we assume that both $m(W)$ and $\varepsilon$ are square integrable random variables. Our goal is to estimate the function $m$ using a dataset $(Y_1,W_1),\ldots,(Y_n,W_n)$ of independent realizations of $(Y,W)$.

Since this framework is a specific case of the more general functional regression framework, usual approaches (which mainly consist in extending classical local methods such as $k$-nearest neighbors, kernel smoothing or local polynomial smoothing) could be used. However, in our context, these methods are known to lead to slow rates of convergence over classical models (see below for detailed references).
Taking advantage of the probabilistic properties of the Wiener coprocess, we aim at defining a new family of models as well as dedicated estimation procedures with faster rates of convergence (in both minimax and adaptive minimax senses). {Despite the fact that considering Brownian paths covariates seems restrictive for pratical purposes, several Brownian diffusion paths could also be considered. Albeilt the systematic theoretical study of such models is beyond the scope of this paper (and is left to further developments), we propose some extensions of our framework as well as some examples of usual processes that can be considered, such as geometric Brownian motions or Ornstein-Uhlenbeck processes.}

In usual functional approaches the set $\mathcal{C}$ is endowed with a metric $d$ \citep[see for example][]{MR2229687,MR2396496,MR2654492} which allows to extend several nonparametric estimators. For example a simple version of the Nadaraya-Watson estimator is given, for any function $w\in\cC$ and any bandwidth $h>0$, by:
\[
\tilde m_h(w) = \sum_{i=1}^n Y_i \frac{\1_{\{d(W_i,w)\leq h\}}}{\sum_{j=1}^n \1_{\{d(W_j,w)\leq h\}}},
\]
where $\1$ stands for the indicator function.
The properties of these estimators are related to the behavior of a quantity known as the \emph{small ball probability} defined for $w\in\cC$ and $h>0$ by $\varphi_w(h)=\p(d(W,w)\le h)$. Pointwise risks of such methods can be generally bounded, up to a positive factor by
\[
h^{\beta} + \left(\frac1{n\varphi_w(h)}\right)^{1/2}
\]
where $\beta$ denotes the smoothness of the mapping $m$ measured in a Hölder sense. For example, if $0<\beta\leq1$ it is assumed that there exists $L>0$ such that $|m(w)-m(w')|\leq Ld(w,w')^\beta$ for any $w,w'\in\mathcal{C}$. Under additional assumptions similar results can be obtained for integrated risks.

The classical assumption $\varphi_w(h) \asymp h^k$ corresponds roughly to the situation where the covariate $W$ lies in some space of finite dimension $k$  \citep[see][]{MR3038004}. This framework corresponds  to the usual nonparametric case. The minimax rates of convergence are then given by $n^{-\beta/(2\beta+k)}$ \citep[see][]{MR2724359}. However if $W$ lies in a functional space, the behavior of $\varphi_w(h)$ is quite different. In our context, where $W$ is a standard Wiener process, it is well-known \citep[see][]{MR1861734} that
\[
\log \p\left(\sup_{t\in[0,1]}|W(t)|\leq h\right) = \log \varphi_0(h) \underset{h\to 0}{\asymp}-h^{-2}
\]
which leads to slower rates of convergence of the form $(\log n)^{-\beta/2}$ assuming a $\beta-$H\"older condition on $m$. We refer the reader to \citet{MR3477653} for recent results with different behavior of $\varphi_w(h)$.

In practical situations, since $\beta$ is unknown, finding adaptive procedures to select the smoothing parameter $h$ is of prime interest. To our best knowledge few papers deal with this problem.
Adaptive procedures based on cross validation have been used in \citet{MR2323791}. \Citet{MR3477653} also propose an adaptation of the method developed by Goldenshluger and Lepski \citep[see][]{GL2011} using an empirical version of the quantity $\varphi_w(h)$. Lower bounds have been investigated by \citet{MR2988463}. In all these papers, the pointwise risk is studied in terms of $\varphi_w(h)$ and theoretical properties are obtained assuming a $\beta$-H\"older condition on $m$ with respect to the metric $d$ with smoothness $\beta\in(0,1]$.

In this paper we follow a different strategy. Taking advantage of probabilistic properties of the Wiener process, similarly to the methodology developed by~\citet{MR3412645} and~\citet{MR3716123},  we consider the Wiener-Itô chaotic decomposition of $m(W)$. Indeed, every random variable that belongs to $\Leb^2_W = \{\km(W) \mid
\km : \cC\to\R \text{ and } \e(\km(W))^2 < +\infty\}$ can be decomposed as a sum of multiple stochastic integrals \citep[see][for more details]{MR2460554}. There exists a unique sequence of functions $(f_\ell)_{\ell\geq 1}$ such that
\begin{equation}\label{eq:WI-chaos}
  m(W) \stackrel{\Leb^2}{=} \e(Y) + \sum_{\ell=1}^\infty \frac1{\ell!} I_{\ell}(f_\ell)(W),
\end{equation}
where $f_\ell$ belongs to $\Lds(\Delta_\ell)$, the set of symmetric and square integrable real-valued functions defined on $\Delta_\ell=[0,1]^\ell$ and
\[
I_\ell(f_\ell)(W) = \int_{\Delta_\ell} f_\ell \,\D W^{\otimes \ell}
=  \int_{\Delta_\ell} f_\ell(u_1,\dotsc,u_\ell) W(\D u_1)\cdots W(\D u_\ell).
\]
{We recall that $f$ is symmetric on $\Delta_\ell$ if for any $(t_1,\dots,t_\ell)\in\Delta_\ell$ and any permutation $\sigma$ of $\{1,2,\dots,\ell\}$, $f(t_1,\dots,t_\ell)=f(t_{\sigma(1)},\dots,t_{\sigma(\ell)})$. } Note that the symmetry implies that the functions $f_\ell$ are isotropic.  
The iterated integral $I_\ell(f_\ell)(W)$ is called a chaos of order $\ell$.

Our approach consists in defining kernel-type estimators $\widetilde{f}_\ell$ of $f_\ell$ using the Itô's isometry, see \eqref{eq:ito-isometry}. Then, based on \eqref{eq:WI-chaos}, we propose the following estimator of $m$
\begin{equation*}
  \hat m_{\mathcal{L}}(W) = \frac1n\sum_{i=1}^n Y_i  + \sum_{\ell=1}^{\mathcal{L}} \frac1{\ell!}I_\ell(\widetilde f_\ell)(W),
\end{equation*}
with $\mathcal L\in\N$. To study these estimators, we  assume that $m$ belongs to a specific class of mappings that satisfy
\[
\sum_{\ell=1}^\infty \frac{e^{2\gamma \ell}}{\ell!} \|f_\ell\|^2_{\Delta_\ell}   \le M^2,
\]
for some $\gamma>0$ and $M>0$, {$\|\cdot\|_{\Delta_\ell}$ is the classical $L_2$ norm on $\Delta_\ell$} and that the $f_\ell$ defined by~\eqref{eq:WI-chaos} are Hölderian. Such classes are quite natural in our context and are connected with the usual Meyer-Watanabe test function space (see section~\ref{sec:scale-mappings} for more details).  

In this case, we find rates of convergence for the prediction error in $\mathbb{L}^p$ norm. Contrary to the classical functional framework, where logarithmic rates are derived, the rates we obtain are intermediate between logarithmic and polynomial rates.

If we assume moreover that the summation in~\eqref{eq:WI-chaos} stops at a known index $L$, we prove that the estimators  $\hat m_{{L}}$ achieve optimal rates of convergence. We derive minimax rates of convergence which are polynomial in $n$ with an exponent that depends on the smoothness of the functions $f_\ell$.  A data-driven procedure, based on the method developed by \cite{GL2011}, is then defined to tune the bandwidths used in the estimation of the functions $f_\ell$. The resulting estimator of $m$ satisfies an oracle-type inequality that allows us to derive adaptive results.

The paper is organized as follows. Section~\ref{sec:2} presents the model and the studied problem. Section~\ref{sec:3} describes the construction of the estimators. Section~\ref{sec:4} gives the main results and Section~\ref{sec:5} is dedicated to the proofs.

\section{Statistical framework}\label{sec:2}

\subsection{Model}\label{sec:scale-mappings}

Let  $W=(W(t) : 0\leq t \leq 1)$ be a standard Brownian motion and let $\varepsilon$ be a centered real-valued random variable independent of $W$. We define:
\[
  Y = m(W) + \varepsilon
\]
where $m:\cC \to \R$ is a given mapping. We assume that $m(W)$ as well as $\varepsilon$ belong to $\Leb^2$, the set of square integrable random variables, then
\begin{equation}\label{eq:WI-chaos2}
  m(W) = \e(Y) + \sum_{\ell=1}^{\infty} \frac1{\ell!} I_{\ell}(f_\ell)(W),
\end{equation}
where for $\ell\in\N$, $f_\ell$ belongs to $\Lds(\Delta_\ell)$.
As mentioned in the introduction we also assume that $f_\ell$ is a regular function. Below we define precisely the functional classes used to measure the smoothness of each function $f_\ell$.

\begin{definition}
    Set $\ell\in\N$ $s_\ell>0$ and  $\Lambda_\ell>0$. The Hölder ball $\cH_\ell(s_\ell, \Lambda_\ell)$  is the set of all functions $f:\Delta_\ell\to\R$ that satisfy the following properties:
    \begin{enumerate}
        \item For any $\alpha=(\alpha_1,\dotsc,\alpha_\ell)\in\N^\ell$ such that $|\alpha|=\sum_i \alpha_i \leq \floor{s_\ell}=\max\{k\in\N \mid k<s_\ell\}$, the partial derivative
        $D^\alpha f$ exists where
        \[
            D^\alpha f = \frac{\partial^{|\alpha|} f}{\partial x_1^{\alpha_1}\dotsc \partial x_\ell^{\alpha_\ell}}.
        \]
        \item For any $x$ and $y$ in $\Delta_\ell$ we have:
        \[
            \sum_{|\alpha|=\floor{s_\ell}}
            \left| D^{\alpha}f(x)-D^{\alpha}f(y) \right|
            \leq \Lambda_\ell |x-y|^{s_\ell-\floor{s_\ell}},
        \]
        where $|\cdot|$ stands for the Euclidean norm of $\R^\ell$.
      \item We have $f\in\Lds(\Delta_\ell)$.
    \end{enumerate}
\end{definition}

Equipped with these notations we can define  a scale of classes $\kA(s,\Lambda,\gamma,M)$ for the mapping~$m$. Roughly, we impose some restrictions on the functions $f_\ell$ that appear in~\eqref{eq:WI-chaos} of two kinds: a minimal smoothness, for each $f_\ell$, is imposed and the growth of the $\mathbf{L}^2$-norm  of the $f_\ell$ is controlled.

\begin{definition}
Set $s=(s_1,s_2,\dotsc)\in(0,+\infty)^\N$, $\Lambda=(\Lambda_1,\Lambda_2,\dotsc)\in(0,+\infty)^\N$, $\gamma>0$ and $M>0$.  We say that $m:\cC\to\R$ belongs to the mapping class $\kA(s,\Lambda,\gamma,M)$  if there exist $a\in\R$ and a sequence of functions $(f_\ell)_{\ell\in\N}$ satisfying 
\begin{equation}\label{eq:WI-chaos2}
  m(W) \stackrel{\Leb^2}{=} a + \sum_{\ell=1}^{\infty} \frac1{\ell!} I_{\ell}(f_\ell)(W),
\end{equation}
with $f_\ell\in\cH_\ell(s_\ell, \Lambda_\ell)$ and
\begin{equation}\label{eq: growth-analytic}
\sum_{\ell=1}^\infty \frac{e^{2\gamma \ell}}{\ell!} \|f_\ell\|^2_{\Delta_\ell}   \le M^2
\end{equation}
where $\|f_\ell\|_{\Delta_\ell} = \left(\int_{\Delta_\ell} f^2_\ell(u) \,\D u\right)^{1/2}$.
\end{definition}

\begin{remark}
  Equation~\eqref{eq: growth-analytic} implies that 
  \begin{equation}
  m(W)\in\bigcap_{\substack{k\geq0\\1< p<e^\gamma+1}} \mathbb{D}_{k,p}
  \end{equation}
  where $\mathbb{D}_{k,p}$ denotes the usual \emph{Sobolev space over the Wiener space} defined in \citet{MR742628}. Note also that, if, for any $\ell\geq1$ we have $\|f_\ell\|_{\Delta_\ell}\leq C^\ell$ for some positive constant $C$, then~\eqref{eq: growth-analytic} is fulfilled for any $\gamma\geq0$.
\end{remark}

{We also define subclasses of the classes $\kA(s,\Lambda,\gamma,M)$ assuming that the summation in \eqref{eq:WI-chaos} stops at a finite index $L\in\N$. }
\begin{definition}
 Set $L\in\N$.
  Set $s=(s_1,\dotsc,s_L)\in(0,+\infty)^L$ and $\Lambda=(\Lambda_1,\dotsc, \Lambda_L)\in(0,+\infty)^L$.  We say that $m:\cC\to\R$ belongs to the mapping class $\kM(s,\Lambda,L,M)$  if there exist $a\in\R$ and a sequence of functions $(f_\ell)_{1\le \ell\le L}$ satisfying 
  \begin{equation}\label{eq:WI-chaos2}
    m(W) = a + \sum_{\ell=1}^{L} \frac1{\ell!} I_{\ell}(f_\ell)(W),
  \end{equation}
  with $f_\ell\in\cH_\ell(s_\ell, \Lambda_\ell)$ and $ \|f_\ell\|^2_{\Delta_\ell} \le M^2 \ell!$ for any $\ell=1,\dotsc,L$.
  \end{definition} 

{More precisely, for any $s$, $\Lambda$ and $M$, the subclasses $\kM(s,\Lambda,L,M)$ satisfy $\kM(s,\Lambda,L,M)\subset \kA(s,\Lambda,\gamma,M')$ with $M'^2=M^2\sum_{\ell=1}^L  e^{-2\gamma\ell}$.} Let us comment on the above definitions since the framework we consider in this paper is quite different to the usual functional framework recalled in the introduction.  In our framework the ``regularity'' of a map $m$ is  seen through the prism of the chaotic decomposition of $m(W)$ and, thus, the functions $f_\ell$. This is not directly linked with the regularity of the mapping $m$ between the space $\mathcal{C}$ endowed with the topology induced by the $\mathbf{L}^2$ norm and $\R$. For example, it can be easily seen that the mapping $m$ defined, for any $w\in\mathcal{C}$ by $m(w) = w^2(1)$ is not continuous (which implies that this function is not hölderian and, thus, cannot be considered in the usual framework). However it is well known that $m(W) = 1 + I_2(1)(W)$. As a consequence, the mapping $m$ falls within our scope since $m$ belongs to $\mathcal{M}(s,\Lambda, 2, 1)$ for any $s$ and $\Lambda$. 

\subsection{{Minimax and adaptive framework}}

The observations consist in a $n$-sample $(Y_1,W_1), \dotsc, (Y_n,W_n)$ distributed as and independent of $(Y,W)$.  Our first goal is to investigate the estimation of $m$, based on these observations, over the classes $\kM(s,\Lambda,L,M)$ and  $ \kA(s,\Lambda,\gamma,M)$ for $0<s\le s^*$ where $\smax$ is fixed. To measure the accuracy of an arbitrary estimator $\tilde m_n = \tilde m(\;\cdot\;; (Y_1,W_1), \dotsc, (Y_n,W_n))$ of $m$, we consider the prediction risk:
\[
  R_p(\tilde m_n, m) = \big(\e|\tilde m_n(W) - m(W)|^p\big)^{1/p}
\]
where $p\geq 2$. The maximal risk of an arbitrary estimator $\tilde m_n$ over a given class of mappings $\kM$ is defined by:
\[
  R_p(\tilde m_n, \kM) = \sup_{\km\in\kM} R_p(\tilde m_n, \km),
\]
whereas the minimax risk is defined, taking the infimum over all possible estimators, by:
\[
  \Phi_n(\kM,p) = \inf_{\tilde m_n} R_p(\tilde m_n, \kM).
\]
An estimator $\tilde m_n$ whose maximal risk is asymptotically bounded, up to a multiplicative factor, by $ \Phi_n(\kM,p)$ is called minimax over $\kM$. Such an estimator is well-adapted to the estimation over $\kM$ but it can perform poorly over another class of mappings. The problem of adaptive estimation consists in finding a single estimation procedure that is simultaneously minimax over a scale of mapping classes. 

Our second goal is to investigate the adaptive estimation of $m$ over the scale of classes 
$\mathcal{M}(s^{*},L,M) = \{\kM_{s,\Lambda,L,M} \mid s\in(0,\smax)^L, \Lambda\in(0,+\infty)^L\}$ where $\smax>0$, $L\in\N$ and $M>0$ are fixed {and known by the statistician}. More precisely our goal is to construct a single estimation procedure $m_n^*$ such that, for any $\kM\in\mathcal M(s^{*},L,M)$, the risk  $R_p(m_n^*, \kM)$ is asymptotically bounded, up to a multiplicative constant, by $\Phi_n(\kM,p)$.
One of the main tools to prove such a result is to find an oracle-type inequality that guarantees that this procedure performs almost as well as the best estimator in a rich family of estimators. Ideally, we would like to have, for any $\km\in\bigcup \kM_{s,\Lambda,L, M}$, an inequality of the following form:
\begin{equation}\label{eq:ideal-oracle}
R_p(m_n^*, \km) \leq \inf_{\eta\in H} R_p(\tilde m_{n,\eta}, \km),
\end{equation}
where  $\{\tilde m_{n,\eta} \mid \eta\in H\}$ is a family of estimators \emph{well-adapted} to our problem in the following sense: for any $\kM\in\mathcal M(s^{*},L,M)$, there exists $\eta\in H$ such that $\tilde m_{n,\eta}$ is minimax over $\kM$.
However, in many situations, \eqref{eq:ideal-oracle} is relaxed and we prove a {weaker} inequality of the type:
\begin{equation}\label{eq:oracle}
R_p(m_n^*, \km) \leq \Upsilon_{\!1,p}\inf_{\eta\in H} R_p^*(\km, \eta) + \Upsilon_{\!2,p} \left(\frac{\log n}{n}\right)^{1/2},
\end{equation}
where $\Upsilon_{\!1,p}$ and $\Upsilon_{\!2,p}$ are two positive constants and  $R_p^*(\km, \eta)$ is an appropriate quantity to be determined that can be viewed as a tight upper bound on $R_p(\tilde m_{n,\eta}, \km)$. Inequalities of the form~\eqref{eq:oracle} are called \emph{oracle-type} inequalities.

Theorems~\ref{theo:2} and~\ref{theo:3} below correspond respectively to an oracle-type inequality and an adaptive result of these types.

\subsection{Extensions to our model}

In this paper, we focus on pure Brownian coprocesses. However our framework allows us to consider a larger class of covariates. Assume that we aim at estimating the regression function $g : \mathcal{C}\to\R$ in the model:
\begin{equation}\label{eq:SDE}
  Y = g(X) + \varepsilon,
\end{equation}
where $X$ is a process driven by the SDE:
\begin{equation}\label{eq:processX}
  \mathrm{d}X_t = b(t, X_t)\mathrm{d}t + \sigma(t, X_t)\mathrm{d}W_t,
  \quad 0\leq t\leq 1.
\end{equation}
Here $\sigma$ and $b$ are assumed to be \emph{known} functions and we also assume that  assumptions guaranteeing the existence and uniqueness of the solution of~\eqref{eq:SDE} are fulfilled. If  for any $0\leq t \leq 1$,  $\sigma(t, X_t)>0$, then, under mild integrability conditions, we have:
\begin{equation}
  \label{eq:reconstructionW}
  W_t =  \int_0^t \frac{\mathrm{d} X_s}{\sigma(s, X_s)}  -  \int_0^t \frac{\mu(s, X_s)}{\sigma(s, X_s)} \mathrm{d} s.
\end{equation}
This implies that there exists a known invertible function $\phi:\mathcal{C}\to\mathcal{C}$ such that  $W = \phi(X)$. In general, this function can be computed by numerical integration. However, in some situations, an exact expression can be obtained using Itô's formula. This is the case for two parametric families of processes widely used to model several practical situations.
First, Ornstein–Uhlenbeck processes $X_t$ are driven by the following SDE:
\begin{equation*}
    \begin{cases}
        \D X_t = -\theta(X_t - \mu) \D t + \sigma \D W_t \\
        X_0 = x_0,
    \end{cases}
\end{equation*}
where $x_0\in\R$ is fixed and $\theta>0$, $\mu\in\mathbb{R}$ and $\sigma>0$ are known parameters. By Itô's formula we have:
\begin{equation*}
    W_t = \frac{1}{\sigma} \left[ X_t - x_0 + \theta\int_0^t \left( X_s - \mu \right) \D s\right].
\end{equation*}
Next, Geometric Brownian motions are used to model stock prices in the Black–Scholes model. Let $\mu\in\R$ and $\sigma>0$ be given parameters. We assume that the process $X=(X_t : t\in I)$ is driven by the following SDE:
\begin{equation}\label{eq:GBS}
    \begin{cases}
        \D X_t = X_t\left( \mu \D t + \sigma \D W_t \right)\\
        X_0 = x_0.
    \end{cases}
\end{equation}
By Itô's formula we have:
\begin{equation*}
    W_t = \frac{1}{\sigma}\left[ \log(X_t/x_0) + 
    \left(\frac{\sigma^2}{2} - \mu\right) t\right].
\end{equation*} 

\begin{remark}
  In practical situation the parameters $\mu$, $\sigma$ and $\theta$ in the above examples are not known. However, estimators of these parameters could be used to estimate the coprocess $W$. This leads to new models where the covariate in observed with errors. The study of such models is beyond the scope of this paper and left to further developments.   
\end{remark}

In view of~\eqref{eq:reconstructionW}, equation~\eqref{eq:SDE} can be written as:
\begin{equation*}
  Y = m(W) + \varepsilon
  \quad\text{where}\quad
  m = g\circ \phi^{-1}.
\end{equation*}
Thus, the regression problem~\eqref{eq:SDE} falls into our framework. The estimation strategy consists of estimating the function $m$ based on the reconstruction of Brownian path $W=\phi(X)$. This can be summarized by the formula:
\begin{equation*}
  \hat g = \widehat{g\circ \phi^{-1}} \circ \phi = \hat m \circ \phi.
\end{equation*} 
Remark also that, in this context, it is relevant to assume that the chaotic decomposition of $m(W)$ is finite. Indeed, under mild assumptions on $b$ and $\sigma$ \citep[see][for more details]{MR1481650}, if $g(X)$ is a polynomial of the terminal value $X_1$ of the process $X$, then the mapping  $m(W) = g(X )$ can be written as a finite chaotic decomposition with smooth functions $f_\ell$.

\section{Estimator construction}\label{sec:3}

In this section we present our estimation procedure. To do so, we first recall classical properties satisfied by Wiener chaos which allow us to construct a family of ``simple'' estimators that depends on a multivariate tuning parameter. Next we construct a procedure which selects, in a data-driven way, this tuning parameter using the methodology developed by \citet{GL2011}.

\subsection{Classical properties of the chaos}

Throughout this paper and in the construction of our statistical procedure, we use the following two fundamental properties satisfied by the iterated integrals.

For $\ell,\ell'\in\mathbb{N}$, Itô's isometry \citep{MR2460554} ensures that, if $g\in\Lds(\Delta_\ell)$ and $g'\in\Lds(\Delta_{\ell'})$, then
\begin{equation}\label{eq:ito-isometry}
    \e\big(I_\ell(g)(W)I_{\ell'}(g')(W)\big) = \delta_{\ell,\ell'} \ell! \int_{\Delta_\ell} g(u)g'(u) \,\D u
\end{equation}
where $\delta_{\ell,\ell'}$ denotes the Kronecker delta.

The hypercontractivity property \citep{MR2962301} will be used to control the concentration of our estimators. Set $q\geq 2$ and $\ell\in\N$. For any $g\in\Lds(\Delta_\ell)$ we have:
\begin{equation}\label{eq:eqcontract}
\left(\e |I_\ell(g)(W)|^q\right)^{1/q}
\leq \kc_\ell(q)
\left(\e I_\ell^2(g)(W)\right)^{1/2}
\quad\text{where}\quad
\kc_\ell(q)=(q-1)^{\ell/2}.
\end{equation}
\subsection{A simple family of estimators}\label{sec:estimators}

Let $\kkk:\R\to\R$ be a function that satisfies the following properties: $\kkk$ is continuous inside $[0,1]$, $\kkk(x)=0$ for any $x\notin[0,1]$,
\[
    \int_0^1 \kkk(x) \,\D x = 1
    \qquad\text{and}\qquad
    \int_0^1 x^s\kkk(x) \,\D x = 0,\quad s=1,\dotsc,\floor{\smax}.
\]
Let $\ell\in\N$. A natural estimator of the function $f_\ell$ is given, for $h\in(0,1)$, by:
\begin{equation}\label{eq:hat-f-ell}
  \hat f^{(\ell)}_{h}(t) = \frac1n \sum_{i=1}^n Y_i I_\ell\left(\Klh{\ell}{h}(t,\cdot)\right)(W_i),
  \qquad t\in\Delta_\ell
\end{equation}
{where} 
$\Klh{\ell}{h}$ is a multivariate kernel defined by:
\[
  \Klh{\ell}{h}(t, u ) = \frac{1}{h^\ell}\prod_{k=1}^{\ell} \kkk\left(\varsigma(t_k)\frac{t_k-u_k}{h}\right)
  \qquad\text{with}\qquad
  \varsigma(\cdot) = 2I_{(1/2,1)} (\cdot) - 1.
\]
This specific construction allows one to obtain an estimator free of boundary bias \citep[see][for more details]{BEEK}.

Indeed, note that for any $t\in\Delta_\ell$ and under regularity assumptions on $f_\ell$ we have:
\begin{align}
f_\ell(t)
  &\stackrel{(h\to0)}{\approx} \int_{\Delta_\ell} f_\ell(u) \Klh{\ell}{h}(t,u)\,\D u\\
  &=
  \e\left(\frac1{\ell!}I_\ell(f_\ell)(W) I_\ell\left(\Klh{\ell}{h}(t,\cdot)\right)(W)\right)\\
  &= \e\left(m(W) I_\ell\left(\Klh{\ell}{h}(t,\cdot)\right)(W)\right)
\end{align}
where the last two lines are obtained using~\eqref{eq:WI-chaos} and~\eqref{eq:ito-isometry}. Since $\varepsilon$ is centered and independent of $W$ we have:
\begin{align}
f_\ell(t)
  &\stackrel{(h\to0)}{\approx}
  \e\left(Y I_\ell\left(\Klh{\ell}{h}(t,\cdot)\right)(W)\right)\\
  &\stackrel{(n\to+\infty)}{\approx} \hat f^{(\ell)}_{h}(t).
\end{align}
Equipped with these notations we define a family of plugin estimators of the mapping $m$. For $\mathcal L\in\N$ and all $\boldsymbol{h}=(h_1,\dotsc,h_{\mathcal{L}})\in(0,1)^{\mathcal{L}}$ we set:
\begin{equation}\label{eq:hat-m-h}
  \hat m_{\boldsymbol{h},\mathcal{L}}(W) = \overline{Y}_{\!\!n} + \sum_{\ell=1}^{\mathcal{L}} \frac{1}{\ell!} I_\ell\big(\hat f^{(\ell)}_{h_\ell}\big)(W)
\end{equation}
where $\overline{Y}_{\!\!n}=\sum_{i=1}^n Y_i/n$.
In the following, we study the rate of convergence of the estimator \eqref{eq:hat-m-h} when $\mathcal{L}=L_n$ where $(L_n)_{n\in\N}$ is a sequence of integers that tends to $\infty$ as $n$ tends to $\infty$ (see Theorem~\ref{thm:minimax}) and $\mathcal{L}=L$ where $L$ is a known fixed integer (see Theorem~\ref{thm:minimax-finite}).

\subsection{Selection procedure}

Set $\smin>0$, $L\in\N$ and $M>0$. Assume that $\mu_4=\left(\e|\varepsilon|^4\right)^{1/4}$ exists. Let $\ell\in\{1,\ldots,L\}$ be fixed and define
$$
\fH_\ell=\left\{h\in(0,1): n^{-1/(2\smin+\ell)}\le h\le (\log n)^{-1} \right\}\cap \{e^{-k}: k\in\N\}.
$$
Now, define
\begin{equation}\label{eq:defMell}
    M(\ell,h) = \frac{\nu(\ell)\left(1+4\sqrt{\log (1/h^\ell)}\right)}{\sqrt{nh^\ell}}
\end{equation}
where
\[
\nu(\ell)=\left(\mu_4+\sum_{k=1}^L \kc_4(k)M \right) \frac{\sqrt{b_{\ell,2}}}{2}
\qquad\text{and}\qquad
\quad b_{\ell,2}=\kc_\ell^2(4)2^\ell \ell!\|\kkk\|_{[0,1]}^{2\ell},
\]
{where the constants $\kc_\ell(k)$ are defined in \eqref{eq:eqcontract}.}
Define for $h\in\fH_\ell$
\[
    B(\ell,h)=\max_{h'\in\fH_\ell}\left\{ \|\hat{f}_{h'}^{(\ell)}-\hat{f}_{h\vee h'}^{(\ell)}\|_{\Delta_\ell} -M(\ell,h') -M(\ell, h\vee h')  \right\}_+
\]
and set
\begin{equation}\label{eq:empirical-tradeoff}
    \hat{h}_\ell= \arg\min_{h\in\fH_\ell} \left\{ B(\ell,h)+M(\ell,h)\right\}.
\end{equation}
The estimation procedure is the defined by $\hat{m}_L=\hat{m}_{\boldsymbol{\hat{h}},L}$ where $\boldsymbol{\hat{h}}=\left(\hat{h}_1,\ldots,\hat{h}_L\right)$.

\begin{remark}
This selection rule follows the principles and the ideas developed by Goldenshluger and Lepski in a series of papers \citep[see][among others]{GL2011,MR3230001}. The quantity $M(\ell, h)$, which is called a \emph{majorant} in the papers cited above, is a penalized version of the standard deviation of the estimator $\hat{f}^{(\ell)}_h$ while the quantity $B(\ell,h)$ is, in some sense, closed to its bias term, see~\eqref{eq:controlBell2}. Finding tight \emph{majorants} is the key point of the method since $\hat{h}_\ell$ is chosen in \eqref{eq:empirical-tradeoff} in order to realize an empirical trade-off between these two quantities.

It is worth noting that the procedure depends on a hyperparameter $\smin>0$ which can be chosen arbitrary small. The introduction of this parameter is due to technical reasons, see~\eqref{eq:controleta1ter} in the proof of Lemma~\ref{lem:5}. This additional assumption (we would like to take $\smin=0$) implies some restrictions on Theorem~\ref{theo:3} below.
\end{remark}

\section{Main results}\label{sec:4}

\subsection{Result for the infinite chaos model}

Our first result studies the risk of our family of estimators over the class $\kA(s,\Lambda,\gamma,M)$. In this class, the function $m$ is decomposed into an infinite sum of chaos:
 \begin{equation}
m(W) \stackrel{\Leb^2}{=} \e(Y) + \sum_{\ell=1}^\infty \frac1{\ell!} I_{\ell}(f_\ell)(W).
\end{equation}

\begin{theorem}\label{thm:minimax}
Set $p\geq 2$, { $\Lambda_*>0$} and $\smax>0$. Set  $s\in(0,\smax)^\N$, { $\Lambda\in(\Lambda_*,+\infty)^\N$} and  $M>0$ and let $\gamma$ be such  { that $2\gamma>\max(2s^* + \log(p-1), \log(3))$}. Assume that {$\mu_p=\e|\varepsilon|^p<+\infty$. }
Define

\begin{equation*}
L_n = \left[( \log n)^{1/2}\right],\qquad C_2=(6(p-1))^{1/2}\|\kkk\|_{[0,1]}
\end{equation*}
where $[\cdot]$ denotes the integer part and  $\boldsymbol{h}_n=\big(h_n^{(\ell)}(s,\Lambda)\big)_{\ell=1,\dotsc,L_n}\in(0,1)^{L_n}$ where for any $\ell=1,\dotsc,L_n$:
\begin{equation*}
h_n^{(\ell)}(s,\Lambda) = \left( \frac{C_2^{2L_n}}{\Lambda_\ell^{2}n} \right)^{\frac{1}{2s_\ell+\ell}}.
\end{equation*}  
There exists a positive constant $\kappa$ depending on $\smax$,  {$\Lambda_*$, $\mu_2$, $\mu_p$,} $\gamma$ and $M$ such that
\begin{equation*}
R_p\big(\hat m_{\boldsymbol{h}_n,L_n},    \kA(s,\Lambda,\gamma,M)    \big)
  \le  \kappa \sum_{\ell=1}^{L_n} \Lambda_\ell^{2\ell/(2s_\ell+\ell)}\left( \frac{C_2^{2L_n}}{n} \right)^{\frac{s_\ell}{2s_\ell+\ell}}.
\end{equation*}
\end{theorem}

Let us briefly comment on this result. Assume first that the parameters $s_\ell$ are constant and denote by $s_0$ their common value. In this case we obtain
\begin{equation*}
  R_p\big(\hat m_{\boldsymbol{h}_n,L_n},    \kA(s,\Lambda,\gamma,M)  \big)
  \le \kappa {(\max_{\ell}\Lambda_\ell^2)}L_n \left( \frac{C_2^{2L_n}}{n} \right)^{\frac{s_0}{2s_0+L_n}}.
\end{equation*}
{This implies that, for $n$ large enough, $R_p\big(\hat m_{\boldsymbol{h}_n,L_n},    \kA(s,\Lambda,\gamma,M)  \big)$, is upperbounded, up to a multiplicative constant by
\begin{equation*}
u_n=
  (\log n)^{1/2} \exp\left(-\frac{s_0}2 (\log n)^{1/2}\right).
\end{equation*}}
Remark that such a rate of convergence lies in-between polylogarithmic rates of convergence and polynomial ones. This result can be compared with those obtained by \citet{MR3412645}. Recall that, in this paper, the authors study a similar model with a Poisson point process covariate. The rates $v_n$ obtained in this paper are slightly better than ours since they obtain, for some $\alpha\in(0,1)$
\begin{equation*}
  \log(v_n) \sim -\frac{\alpha}2\left( \log n\right)^{1/2} \left(\log\log n \right)^{1/2}
\end{equation*}
whereas, in our case,
\begin{equation*}
  \log(u_n) \sim -\frac{s_0}2\left( \log n\right)^{1/2}.
\end{equation*}
However remark that their study is limited to $p=2$ and $s_0=1$ and that, moreover, they assume that the response $Y$ is a bounded variable. In our situation neither $m(W)$ nor $\varepsilon$ are assumed to be bounded.

\subsection{Results for finite chaos model}

In the three following results, we assume that it exists {a known integer $L\in \N$} such that 
\begin{equation}
m(W) \stackrel{\Leb^2}{=} \e(Y) + \sum_{\ell=1}^L \frac1{\ell!} I_{\ell}(f_\ell)(W).
\end{equation}
Our second result proves that the minimax rate of convergence over the class $\kM_{s,\Lambda,L,M}$ is of the same order as:
\begin{equation}\label{eq:minimax-rate}
  \phi_n(s,\Lambda) = \max\left\{ \Lambda_\ell^{2\ell/(2s_\ell+\ell)}n^{-s_\ell/(2s_\ell+\ell)}\mid \ell=1,\dotsc,L\right\}.
\end{equation}
\begin{theorem}\label{thm:minimax-finite}
  Set $p\geq 2$, $s\in(0,\smax)^L$, $\Lambda\in(\Lambda_*,+\infty)^L$, $M>0$ and assume that $\mu_p = (\e|\varepsilon|^p)^{1/p}<+\infty$. Define  $\tilde{\boldsymbol{h}}_n=\big(\tilde h_n^{(\ell)}(s,\Lambda)\big)_{\ell=1,\dotsc,L}\in(0,1)^L$ where:
  \[
    \tilde h_n^{(\ell)}(s,\Lambda) = \left(\frac{1}{\Lambda_{\ell}^2\,n}\right)^{\frac{1}{2s_\ell+\ell}},
    \qquad \ell=1,\dotsc,L.
  \]
There exist two positive constants $\kappa_*$ and $\kappa^*$ that depend only on $L$, $M$, {$\Lambda_*$, $\mu_2$, $\mu_p$} and $\smax$ such that
  \begin{equation}
    \label{eq:minimax-upper-bound}
    \limsup_{n\to+\infty} \phi_n^{-1}(s,\Lambda) R_p\left(\hat m_{\tilde{\boldsymbol{h}}_n,L}, \kM_{s,\Lambda,L,M}\right) \leq \kappa^*
  \end{equation}
  and
  \begin{equation}\label{eq:minimax-lower-bound}
    \liminf_{n\to+\infty} \phi_n^{-1}(s,\Lambda) \Phi_n\left(\kM_{s,\Lambda,L,M},p\right) \geq\kappa_*.
  \end{equation}
\end{theorem}

Note that this result also ensures that the family of estimators constructed in Section~\ref{sec:estimators} is \emph{well-adapted} to our problem. The next result states an \emph{oracle-type} inequality satisfied by our data-driven estimator $\hat{m}_L$.

\begin{theorem}\label{theo:2} Set $p\geq2$ and assume that 
for any $\ell=1,\ldots,L$, $\|f_\ell\|_{\Delta_\ell}^2\le M^2\ell!$ and that for any $q\ge 1$ the moment $\mu_q=\left(\e|\varepsilon|^q\right)^{1/q}$ exists.
Then:
\[
    R_p(\hat{m}_L, m)\le \Upsilon_1\sum_{\ell=1}^L \inf_{h\in\fH_\ell}  \left[ \max_{\substack{h'\in\fH_\ell\\ h'\le h}} \|\e \hat{f}_{h'}^{\ell}-f_\ell\|_{\Delta_\ell} +M(\ell,h) \right] +\Upsilon_2\left(\frac{\log n}{n}\right)^{1/2},
\]
where $\Upsilon_1$ and $\Upsilon_2$ are two positive constants that depend on $L$, $M$, $\smin$ and $\smax$.
\end{theorem}

Using Theorems~\ref{thm:minimax-finite} and~\ref{theo:2} we can derive our last result: the data-driven estimation procedure is adaptive, up to a logarithmic factor, over the scale $\{\kM_{s,\Lambda, L, M}:s\in (\smin,s^*)^L, \Lambda\in(0,+\infty)^L, M>0\}$.

\begin{theorem}\label{theo:3} Set $p\geq 2$ and assume that for any $q\ge 1$ the moment $\mu_q=\left(\e|\varepsilon|^q\right)^{1/q}$ exists. For any $s\in (\smin,s^*)^L$, any $\Lambda\in(0,+\infty)^L$, any $M>0$, we have
\[
  \limsup_{n\to+\infty} \tilde{\phi}^{-1}_n(s) R_p(\hat{m}_L, \kM_{s,\Lambda, L, M})\le \kappa^{**}
\]
where $\kappa^{**}$ is a positive constant that depends on $\Lambda$, $L$, $M$, $\smin$ and $\smax$ and
$$\tilde{\phi}_n(s)=\max\left\{ \left(\frac{\log n}{n}\right)^{s_\ell/(2s_\ell+\ell)}\mid \ell=1,\dotsc,L\right\}.$$
\end{theorem}

\begin{remark}
While the selection procedure is defined using the $\fL^2$-norms, the procedure is adaptive for any $p\geq 2$. This phenomenon is due to the hypercontractivity property, see~\eqref{eq:controlmhat}. 
Note that in Theorem~\ref{theo:2}, the quantity $$\max_{\substack{h'\in\fH_\ell\\ h'\le h}} \|\e \hat{f}_{h'}^{\ell}-f_\ell\|_{\Delta_\ell}$$ is a tight upper bound of the bias term of the estimator $\hat{f}_{h}^\ell$.

This result ensures that our data-driven procedure is adaptive, up to a logarithmic factor, over a large scale of mapping classes.

The presence of the extra logarithmic factor in the adaptive rate of convergence is not usual for prediction risks. This term is introduced in the definition of $M(\ell, h)$ to control the deviation of the estimator \eqref{eq:hat-f-ell} based on the variables $I_\ell\left(K_h^{(\ell)}(t,\cdot)\right)(W_i)$. See \eqref{eq:explicationlogn} for more details.
\end{remark}

\section{Proofs}\label{sec:5}

We first consider some notations and lemmas. Define for $i\in\{1,\ldots,n\}$, $\ell\in\{1,\ldots,L\}$ and $h\in (0,1)$
\[
    \xi_{i,\ell}(t,h) = I_\ell\big(K_{h}^\ell(t,\cdot)\big)(W_i),
    \quad\text{and}\quad
    \xi_{\ell}(t,h) = I_\ell\big(K_{h}^\ell(t,\cdot)\big)(W)
\]
and
\[
    \Theta_{i,\ell} = I_\ell(f_\ell)(W_i)
    \quad\text{and}\quad
    \Theta_{\ell} = I_\ell(f_\ell)(W).
\]

\begin{lemma}\label{lem:1}
We have, for any $\ell\in\{1,\ldots,L\}$, $h\in\fH_\ell$ and  $r\geq 1$
\[\left(\e \xi_{\ell}^{2r}(t,h)\right)^{1/r} \leq b_{\ell,r} h^{-\ell}
\quad \text{with} \quad b_{\ell,r}=\kc_\ell^2(2r)2^\ell \ell!\|\kkk\|_{[0,1]}^{2\ell}.\]
Moreover for $\varphi>0$ and $q\geq 1$
\[
\e\left(|\xi_{\ell}(t,h)|^r\1_{|\xi_{\ell}(t,h)| > \varphi}\right)
\leq  (b_{\ell,r})^{r/2} (b_{\ell,q})^{q/2}\varphi^{-q}h^{-\ell(r+q)/2}.
\]
\end{lemma}
\begin{lemma}\label{lem:5}
Let $\ell\in\{1,\ldots,L\}$ and $h\in\fH_\ell$. Let $\chi,\chi_1,\ldots,\chi_n$ be i.i.d random variables such that, for any $r\ge 1$
\begin{equation}\label{eq:hypchi}
\left(\e |\chi|^{2r}\right)^{1/(2r)}\le a_r <+\infty.
\end{equation}
Define
\[
    \fU(t)=\frac{1}{n} \sum_{i=1}^n\bigg\{
    \chi_i \xi_{i,\ell}(t,h)
    - \e\left(
    \chi_i \xi_{i,\ell}(t,h)
    \right)
    \bigg\}
\]
and
\[
    T= (1+\delta) \frac{a_2\sqrt{b_{\ell,2}}}{2\sqrt{nh^\ell}}
    \qquad\text{with}\qquad
    \delta = 4\sqrt{\log (1/h^\ell)}.
\]
Then there exists a positive constant $C>0$ such that
\begin{equation}\label{eq:conclusionlemma5}
    \e\left\{\|\fU\|_{\Delta_\ell}-T\right\}_+^2\le Cn^{-1}.
\end{equation}
\end{lemma}

The following Lemma recalls the Bousquet's version of Talagrand's concentration  inequality \citep[see][]{MR1890640,MR3185193}.
\begin{lemma}[Bousquet's inequality]\label{lem:bousquet}
Let $X_1,\ldots, X_n$ be independent identically distributed random variables. Let $\mathcal{S}$ be a countable set of functions and define $Z=\sup_{s\in\mathcal{S}}\sum_{i=1}^n s(X_i)$. Assume that, for all $i=1,\ldots, n$ and $s\in\mathcal{S}$, we have $\e s(X_{i})=0$ and $s(X_{i})\le 1$ almost surely. Assume also that $v=2\e Z+\sup_{s\in\mathcal{S}}\sum_{i=1}^n \e (s(X_i))^2<\infty $. Then we have for all $t>0$
\begin{equation}\label{eq:bousquet}
\p\left(Z-\e Z\ge t\right) \leq \exp\left\{ -\frac{t^2}{2(v+\frac{t}{3})} \right\}.
\end{equation}
\end{lemma}

\subsection{Proof of Theorem~\ref{thm:minimax}}

Set $p\geq 2$, {$s\in(0,\smax)^{\N}$, $\Lambda\in(\Lambda_*,+\infty)^{\N}$}, $\gamma>0$, and  $M>0$.
For the sake of readability we denote $\boldsymbol{h}=\boldsymbol{h}_n$, $h_\ell=h^{(\ell)}_n(s,\Lambda)$, $K_\ell = \Klh{\ell}{h_\ell}$, $\xi_\ell(t)=\xi_\ell(t, h_\ell)$ and $\hat f_\ell = \hat f^{(\ell)}_{h_\ell}$. 

\paragraph{Decomposition of the risk.}
Using the triangle inequality we have:
\begin{align}
    R_p(\hat m_{\boldsymbol{h},L_n},m)
    &\leq
      \big(\e|\bar Y_n-\e(Y)|^p\big)^{1/p}
    + \sum_{\ell=1}^{L_n} \frac{1}{\ell!}\left(\e\left|I_\ell\big(\hat f_{\ell}-f_{\ell}\big)(W)\right|^p\right)^{1/p}\\
    &\qquad+ \sum_{\ell=L_n+1}^{+\infty} \frac{1}{\ell!}
    \left(\e\left|I_\ell\big(f_{\ell}\big)(W)\right|^p\right)^{1/p}\\
    &\leq \big(\e|\bar Y_n-\e(Y)|^p\big)^{1/p}
    + \sum_{\ell=1}^{L_n} \frac{\kc_\ell(p)}{\ell!}\left(\e\left|I_\ell\big(\hat f_{\ell}-f_{\ell}\big)(W)\right|^2\right)^{1/2}\\
    &\qquad+ \sum_{\ell=L_n+1}^{+\infty} \frac{\kc_\ell(p)}{\ell!}
    \left(\e\left|I_\ell\big(f_{\ell}\big)(W)\right|^2\right)^{1/2}
\end{align}
Last line comes from the hypercontractivity property. Now, using Itô's isometry, we obtain:
\begin{align}
  R_p(\hat m_{\boldsymbol{h},L_n},m)
  &\leq \big(\e|\bar Y_n-\e(Y)|^p\big)^{1/p}
    + \sum_{\ell=1}^{L_n} \frac{\kc_\ell(p)}{\sqrt{\ell!}}
    \left(\e \|\hat f_{\ell}-f_{\ell}\|^2_{\Delta_\ell}\right)^{1/2}\\
    &\qquad +\sum_{\ell=L_n+1}^{+\infty} \frac{\kc_\ell(p)}{\sqrt{\ell!}}
    \|f_{\ell}\|_{\Delta_\ell}\\
  &\leq \big(\e|\bar Y_n-\e(Y)|^p\big)^{1/p}+ \sum_{\ell=1}^{L_n} \frac{\kc_\ell(p)}{\sqrt{\ell!}}  \big(B(\ell) + V(\ell)\big)\\
  &\qquad + \sum_{\ell=L_n+1}^{+\infty} \frac{\kc_\ell(p)}{\sqrt{\ell!}}
    \|f_{\ell}\|_{\Delta_\ell}\label{eq:risk-decomposition}
\end{align}
where the \emph{bias term} $B(\ell)$ and the \emph{stochastic term} $V(\ell)$ are defined by:
\[
    B(\ell) = \|\e\hat f_{\ell}-f_{\ell}\|_{\Delta_\ell}
    \quad\text{and}\quad
    V(\ell) =\left( \e\|\hat f_{\ell}-\e\hat f_{\ell}\|^2_{\Delta_\ell}\right)^{1/2}.
\]

\paragraph{Study of the constant term.}
Remark that
\begin{align}
  \left(\e|\bar Y_n-\e Y|^p\right)^{1/p}
  &= \left(\e\left|\frac1n \sum_{i=1}^n (Y_i-\e Y_i)\right|^p\right)^{1/p}\\
  &\le \frac{C_{1,p}}{n^{1-1/p}} \sigma_Y(p) + \frac{C_{2,p}}{n^{1/2}}\sigma_Y(2) \label{eq:bornerosen}
\end{align}
where the last line is obtained using Rosenthal's inequality \citep[]{MR770640}. Here $C_{1,p}$ and $C_{2,p}$ denote two positive constants while $\sigma_Y(p) = (\e|Y - \e Y|^p)^{1/p}$.
Moreover since
\begin{equation}
    Y - \e Y
    = \sum_{\ell = 1}^{+\infty} \frac{\Theta_\ell}{\ell!} + \varepsilon
\end{equation}
the hypercontractivity property, implies that, for any $p\geq 2$
\begin{align}
 \sigma_Y(p)
 &\leq
\sum_{\ell = 1}^{+\infty} \frac1{\ell!}\left(\e\left| \Theta_\ell\right|^p\right)^{1/p}
 + \left(\e|\varepsilon|^p\right)^{1/p}\\
 &\leq \sum_{\ell = 1}^{+\infty} \frac{\kc_\ell(p)}{\ell!}\left(\e \Theta_\ell^2\right)^{1/2}
 + \left(\e|\varepsilon|^p\right)^{1/p}\\
 &= \sum_{\ell = 1}^{+\infty} \frac{\kc_\ell(p)}{\sqrt{\ell!}}\|f_\ell\|_{\Delta_\ell}
 + \left(\e|\varepsilon|^p\right)^{1/p}
\end{align}
Last line comes from Itô's isometry. Now, using that $m\in\kA(s,\Lambda,\gamma,M) $ and applying Cauchy-Schwarz inequality we obtain: 
\begin{equation}
  \sigma_Y(p)
  \leq M \left( \sum_{\ell=1}^{+\infty} \exp(-2\gamma \ell) \kc_\ell^2(p) \right)^{1/2}+(\e|\varepsilon|^p)^{1/p} = s_Y(p)
\end{equation}
where, using the definition of $\kc_\ell(p)$ and the fact that $2\gamma>\log(p-1)$
\[
s_Y(p)
=M \left( \sum_{\ell=1}^{+\infty} \exp(-2\gamma \ell) (p-1)^{\ell} \right)^{1/2} +(\e|\varepsilon|^p)^{1/p}
<+\infty.
\]
We finally obtain
\begin{align}
  \left(\e|\bar Y_n-\e Y|^p\right)^{1/p}
  &\le \frac{C_{1,p}s_Y(p)}{n^{1-1/p}} + \frac{C_{2,p}s_Y(2)}{n^{1/2}}\\
  &\le \frac{\kappa_0}{n^{1/2}},\label{eq:kappa0}
\end{align}
with $\kappa_0 = C_{1,p}s_Y(p)+C_{2,p}s_Y(2)$ depends only on $M$, { $\gamma$}, $\mu_2$ and $\mu_p$.

\paragraph{Study of the bias term.}

Set $\ell\in\{1,\dotsc,L_{n}\}$ and note that:
\[
    \e\hat f_{\ell}(t)
    = \e\left((Y\xi_\ell(t)\big)(W)\right)
    = \e\left(m(W)\xi_\ell(t)\right)
    = \frac{1}{\ell!}\e\left(\Theta_\ell \xi_\ell(t)\right)
\]
Using Itô's isometry we thus obtain:
\[
\e\hat f_{\ell}(t) =  \int_{\Delta_\ell} f_\ell(u) K_\ell(t,u) \,\D u.
\]
To apply multivariate Taylor formula we introduce, for any $\alpha=(\alpha_1,\dotsc,\alpha_ \ell)\in(\N\cup\{0\})^\ell$, the notation $|\alpha|=\alpha_1+\dotsc+\alpha_\ell$. Moreover we define $s_\ell = m_\ell+\gamma_\ell$ with $m_\ell\in\N\cup\{0\}$ and $0 < \gamma_\ell\leq 1$.
Since $f_\ell\in\cH_\ell(s_\ell,\Lambda_\ell)$, we obtain, using classical arguments \citep[see][]{BEEK}, that:
\begin{align}\label{eq:bias-term}
  B(\ell)
  &\leq \left[ 2^\ell m_\ell
  \sum_{|\alpha|=m_\ell}
  \prod_{i=1}^\ell\int_{0}^1 |\kkk(y)| y^{\alpha_i+\gamma_\ell} \,\D y \right] \Lambda_\ell h_\ell^{s_\ell}\nonumber\\
  &\leq \kb_\ell(\kkk,s)\Lambda_\ell h_\ell^{s_\ell}
\end{align}
where
\[
  \kb_\ell(\kkk,s) = \big(\lfloor s^* \rfloor\wp(\lfloor s^* \rfloor)\big)\big(2\|\kkk\|_{[0,1]}\big)^\ell  \ell^{\lfloor s^* \rfloor}
\]
and $\wp(\cdot)$ denotes the partition function of an integer.
We then obtain:
\begin{equation*}
  \sum_{\ell=1}^{L_n} \frac{\kc_\ell(p)}{\sqrt{\ell!}}  B(\ell)
  \leq \big(\lfloor s^* \rfloor\wp(\lfloor s^* \rfloor)\big)\sum_{\ell=1}^{L_n} \frac{\big({2(p-1)^{1/2}\|\kkk\|_{[0,1]}}\big)^{\ell}}{\sqrt{\ell!}}\Lambda_\ell h_\ell^{s_\ell}
\end{equation*}
Since the sequence
\begin{equation*}
  \frac{\big({2(p-1)^{1/2}\|\kkk\|_{[0,1]}}\big)^{\ell}}{\sqrt{\ell!}}
\end{equation*}
tends to $0$ as $\ell$ goes to infinity, there exists an absolute constant $C_0>0$ that depends only on $p$, $s^*$ and $\|\kkk\|_{[0,1]}$ such that:
\begin{equation}\label{eq:bound-bias}
  \sum_{\ell=1}^{L_n} \frac{\kc_\ell(p)}{\sqrt{\ell!}}  B(\ell) 
  \leq C_0 
  \sum_{\ell=1}^{L_n} \Lambda_\ell h_\ell^{s_\ell}.
\end{equation} 

\paragraph{Study of the stochastic term $V(\ell)$}

Set $\ell\in\{1,\dotsc,L_n\}$. We have:
\begin{align}
    V(\ell)
    &= \left(\frac{\fV_1(\ell)+\fV_2(\ell)}n\right)^{1/2}\\
    &\le \frac{(\fV_1(\ell))^{1/2}+(\fV_2(\ell))^{1/2}}{n^{1/2}}
\end{align}
where
\[
  \fV_1(\ell) =  \int_{\Delta_\ell} \e\left(\varepsilon^2 \xi_\ell^2(t)\right) \,\D t
\]
and
\[
  \fV_2(\ell) = \int_{\Delta_\ell} \var\left(m(W) \xi_{\ell}(t) \right) \,\D t.
\]
Then we have 
\begin{equation}
\sum_{\ell=1}^{L_n} \frac{\kc_\ell(p)}{\sqrt{\ell!}}  V(\ell) 
  \leq \sum_{\ell=1}^{L_n} \kc_\ell(p)\left(\frac{V_1(\ell)}{n\ell!}\right)^{1/2}+\sum_{\ell=1}^{L_n} \kc_\ell(p)\left(\frac{V_2(\ell)}{n\ell!}\right)^{1/2}
\end{equation}
Since $W$ and $\varepsilon$ are independent Lemma~\ref{lem:1} implies:
\begin{equation}
\fV_1(\ell)
= \int_{\Delta_\ell} \e\left(\varepsilon^2\right) \e\left(\xi_\ell^2(t)\right) dt\leq \frac{\mu_2^2 b_{\ell,1}}{h_\ell^\ell}.
\end{equation}
Then we have,
\begin{align}\label{eq:V1}
  \sum_{\ell=1}^{L_n} \kc_\ell(p)\left(\frac{V_1(\ell)}{n\ell!}\right)^{1/2}\le \mu_2\sum_{\ell=1}^{L_n} a_1^{\ell}\left(\frac{1}{nh_\ell^\ell}\right)^{1/2}
  \le \mu_2a_1^{L_n} \sum_{\ell=1}^{L_n} \left(\frac{1}{nh_\ell^\ell}\right)^{1/2}
\end{align}
where
\[
a_1=(2(p-1))^{1/2}\|\kkk\|_{[0,1]}\geq 1.  
\]
Now we have:
\begin{align}
    \fV_2(\ell)
    &= \int_{\Delta_\ell}
    \e\left(\sum_{k=1}^{L_n}\frac1{k!}\Theta_k\xi_\ell(t)\right)^2 \,\D t\\
    &=\sum_{k,k'=1}^{L_n}\frac1{k!k'!}
    \int_{\Delta_\ell}
    \e\left(\Theta_k\Theta_{k'} \xi_\ell^2(t)\right) \,\D t.
\end{align}
Using Cauchy-Schwarz inequality we obtain:
\[
    \fV_2(\ell)
    \leq \sum_{k,k'=1}^{L_n}
    \frac{\left(\e \Theta_k^4\right)^{1/4}\left(\e \Theta_{k'}^4\right)^{1/4}}{k!k'!}
    \int_{\Delta_\ell}
    \left(\e \xi_\ell^4(t)\right)^{1/2} \,\D t.
\]
Now, using Lemma~\ref{lem:1} and that $m\in\kA(s,\Lambda,\gamma,M) $, 
we obtain:
\begin{align}
  \fV_2(\ell)
&\leq
\left(\sum_{k}^{L_n}\frac{1}{\sqrt{k!}}\kc_k(4) \|f_k\|_{\Delta_k}\right)^2b_{\ell,2}h_ \ell^{-\ell}\\
& \leq M^2 \left(\sum_{k=1}^{\infty}e^{-2\gamma k}\kc_k^2(4)\right)b_{\ell,2}h_ \ell^{-\ell}.
\end{align}
This implies that 
\begin{align}\label{eq:V2}
  \sum_{\ell=1}^{L_n} \kc_\ell(p)\left(\frac{V_2(\ell)}{n\ell!}\right)^{1/2}\le a_2\sum_{\ell=1}^{L_n} C_2^{\ell}\left(\frac{1}{nh_\ell^\ell}\right)^{1/2}
  \le a_2C_2^{L_n} \sum_{\ell=1}^{L_n} \left(\frac{1}{nh_\ell^\ell}\right)^{1/2}
\end{align}
where
\[
C_2=\sqrt{3}a_1
\qquad\text{and}\qquad
a_2=M  \left(\sum_{k=1}^{\infty}e^{-2\gamma k}\kc_k^2(4)\right)^{1/2}<\infty. 
\]
Note that $a_2$ is finite since $2\gamma>\log(3)$
Combining \eqref{eq:V1} and \eqref{eq:V2}, we obtain, denoting $C_1= a_2+\mu_2$, that 
\begin{equation}\label{eq:bound-var}
  \sum_{\ell=1}^{L_n} \frac{\kc_\ell(p)}{\sqrt{\ell!}}  V(\ell) 
    \leq C_1 C_2^{L_n} \sum_{\ell=1}^{L_n}  \left(\frac{1}{nh_\ell^\ell}\right)^{1/2}. 
\end{equation}

\paragraph{General bound on the risk}

Combining \eqref{eq:kappa0}, \eqref{eq:V1} and \eqref{eq:V2}, the following bound can be easily obtained:
\begin{align*}
  R_p(\hat m_{\mathbf{h},L_n}, m)
  \leq &\kappa_0 n^{-1/2} + C_0 
  \sum_{\ell=1}^{L_n} \Lambda_\ell h_\ell^{s_\ell}+ C_1 C_2^{L_n} \sum_{\ell=1}^{L_n}  \left(\frac{1}{nh_\ell^\ell}\right)^{1/2}\\ 
  &\qquad+\sum_{\ell=L_n+1}^{+\infty} \frac{\kc_\ell(p)}{\sqrt{\ell!}}
  \|f_{\ell}\|_{\Delta_\ell}
\end{align*}

\paragraph{Study of the residual term}

Finally we have using that $m\in\kA(s,\Lambda,\gamma,M) $
\begin{align}
  \sum_{\ell=L_n+1}^{+\infty} \frac{\kc_\ell(p)}{\sqrt{\ell!}}
  \|f_{\ell}\|_{\Delta_\ell}
  &\leq   \psi_\gamma(L_n)
  \left(\sum_{\ell=L_n+1}^{+\infty} (p-1)^{\ell}e^{-2\gamma\ell}\right)^{1/2}\\
  &\leq \frac{\psi_\gamma(L_n)}{\sqrt{2\gamma_p}} e^{-\gamma_p L_n}\label{eq:residual}
\end{align}
where
\begin{equation*}
  \gamma_p = \gamma - \frac{\log(p-1)}2
  \quad\text{and}\quad
  \psi_\gamma(L) = \left( \sum_{\ell=L+1}^{+\infty}  \frac{\|f_{\ell}\|_{\Delta_\ell}^2 e^{2\gamma\ell}}{\ell!}
   \right)^{1/2}. 
\end{equation*}
Note that $\psi_\gamma(L)$ tends to $0$ as $L$ tends to infinity. 

\paragraph{Upper bound}

Using the definitions of $L_n$ and $h_\ell$ we have 
\begin{equation*}
  R_p(\hat m_{\mathbf{h},L_n}, m)
  \leq \kappa_0 n^{-1/2} + (C_0+C_1) \sum_{\ell=1}^{L_n} \Lambda_\ell^{2\ell/(2s_\ell+\ell)}\left( \frac{C_2^{2L_n}}{n} \right)^{\frac{s_\ell}{2s_\ell+\ell}}
  + \frac{\psi_\gamma(L_n)}{\sqrt{2\gamma_p}} e^{-\gamma_p L_n}.
\end{equation*}
{Now, remark that, since $\gamma_p>s^{*}$ and $C_2>1$,
we have 
\begin{align*}
  e^{-\gamma_p L_n}& \le   e^{-s^{*} L_n}\le e^{s^*}e^{-s^*(\log n)^{1/2}}=e^{s^*}e^{-s^*\frac{\log n}{L_n}\left[\frac{L_n+1}{(\log n)^{1/2}}-\frac{1}{(\log n)^{1/2}}\right]}
  \\&\le e^{3s^*}e^{-\frac{s^{*} \log n}{L_n}} \le e^{3s^*} e^{-\frac{s^{*} \log n}{2s^{*}+L_n}}=e^{3s^*}n^{-\frac{s^{*}}{2s^{*}+L_n}}\le e^{3s^*}\left( \frac{C_2^{2L_n}}{n} \right)^{\frac{s^*}{2s^*+L_n}}.
  \end{align*}
Then there exists a positive constant $a_4$ that depends on $\Lambda_*$ and $s^*$ such that}
\begin{equation*}
  e^{-\gamma_p L_n} 
  \leq a_4 \sum_{\ell=1}^{L_n} \Lambda_\ell^{2\ell/(2s_\ell+\ell)}\left( \frac{C_2^{2L_n}}{n} \right)^{\frac{s_\ell}{2s_\ell+\ell}}.
\end{equation*}
This implies that
\begin{equation*}
  R_p(\hat m_{\mathbf{h},L_n}, m) 
  \leq (C_0+C_1) \sum_{\ell=1}^{L_n} \Lambda_\ell^{2\ell/(2s_\ell+\ell)}\left( \frac{C_2^{2L_n}}{n} \right)^{\frac{s_\ell}{2s_\ell+\ell}} + \rho_n
\end{equation*}
where $\rho_n$ is a negligeable reminder term.

\subsection{Proof of Theorem~\ref{thm:minimax-finite}}

This proof is decomposed into two parts. We first prove the upper bound~\eqref{eq:minimax-upper-bound} and then the lower bound~\eqref{eq:minimax-lower-bound}.

\subsubsection{Proof of the upper bound}

For the sake of readability we denote $\tilde{\boldsymbol{h}}=\tilde{\boldsymbol{h}}_n$ and $h_\ell=\tilde h^{(\ell)}_n(s,\Lambda)$. Following the same notations as in the proof of Theorem~\ref{thm:minimax}, we have
\begin{equation}
    R_p(\hat m_{\tilde{\boldsymbol{h}},L},m)
    \leq \big(\e|\bar Y_n-\e(Y)|^p\big)^{1/p}+ \sum_{\ell=1}^{L} \frac{\kc_\ell(p)}{\sqrt{\ell!}}  \big(B(\ell) + V(\ell)\big).
\end{equation}
Note that in this case there is no residual term.
Similarly to the proof of Theorem~\ref{thm:minimax}, and using the same notations, we have
\begin{equation*}
\big(\e|\bar Y_n-\e(Y)|^p\big)^{1/p}\le \kappa_0n^{-1/2},
\end{equation*}
with $\kappa_0$ depending on $M$, $\mu_2$ and $\mu_p$.
The bias term  satisfies
\begin{equation*}
  \sum_{\ell=1}^{L} \frac{\kc_\ell(p)}{\sqrt{\ell!}} B(\ell)
  \leq 
  C_0 \sum_{\ell=1}^L \Lambda_\ell h_\ell^{s_\ell}
\end{equation*} and the stochastic term satisfies
\begin{equation*}
\sum_{\ell=1}^{L} \frac{\kc_\ell(p)}{\sqrt{\ell!}}  V(\ell)\le C_4 \sum_{\ell=1}^L \left(\frac{1}{nh_\ell^\ell}\right)^{1/2},
\end{equation*}
where $C_4$ depends on $L$, $M$, and $\mu_2$.
Now  by  substituting $h_\ell$ by its value, we obtain
\begin{align}
    R_p(\hat m_{\tilde{\boldsymbol{h}},L},m)
    &\leq \kappa^* \max\left\{\Lambda_\ell^{2\ell/(2s_\ell+\ell)} n^{-s_\ell/(2s_\ell+\ell)} \mid \ell=1,\dotsc,L\right\}
\end{align}
where $\kappa^*$ is a positive constant that depends only on $L$, $M$ {$\mu_2$, $\mu_p$, $\Lambda_*$} and $\smax$.
This ends the proof of the upper bound. Now, let us prove the lower bound.

\subsubsection{Proof of the lower bound}

Note that for any $p\geq2$ and any estimator $\tilde m_n$ of $m$ we have $R_p(\tilde m_n, m) \geq R_2(\tilde m_n, m)$. This implies that, to prove the lower bound, it is sufficient to consider the case $p=2$.

\paragraph{Method.}

We fix $s\in(0,\smax)^L$, $\Lambda\in(0,+\infty)^L$ and $M>0$. To prove the lower bound over the space $\kM(s,\Lambda,L,M)$, we define
\[
\ell = \argmax_{k=1,\dotsc,L} \Lambda_k^{2k/(2s_k+k)}n^{-s_k/(2s_k+k)}
\]
and we follow the strategy developed by~\citet{MR3716123}. In particular Lemma~6.1 of this paper implies (using Itô's isometry combined with Theorem~2.5 in~\citet{MR2724359}) that the problem boils down to find a finite family of functions $\left\{g_\omega\right\}_{\omega\in \mathcal{W}}$ with cardinal $|\mathcal{W}|\ge 2$ that satisfies the following assumptions:
\begin{itemize}
\item[(i)] the null function $0\in\{g_\omega\}_{\omega\in\mathcal{W}}$.
    \item[(ii)] for any $\omega\in \mathcal{W}$, the function $g_\omega\in \cH_\ell(s_\ell, \Lambda_\ell, M)$ and $\|g_\omega\|_{\Delta_\ell}^2\leq \ell! M^2$
    \item[(iii)] there exists $\kappa_*>0$ such that for $\omega\neq\omega'$, $\|g_\omega-g_{\omega'}\|_{\Delta_\ell}\ge 2\kappa_*\phi_n(s,\Lambda)$
    \item[(iv)] there exists $0<\alpha<1/8$ such that
    \[
    \frac{n}{|\mathcal{W}|}\sum_{\omega\in\mathcal{W}}\| g_{\omega}\|_{\Delta_\ell}^2\le 2\alpha\log(|\mathcal{W}|).
    \]
\end{itemize}
Under these assumptions, the lower-bound \eqref{eq:minimax-lower-bound} holds for $p=2$.

\paragraph{Notation.} Here, we construct a finite set of functions used in the rest of the proof. We consider the function $\psi:\R\to\R$ defined, for any $u\in\R$ by
\[
	\psi(u)=\exp(-1/(1-u^2))\1_{(-1,1)}(u).
\]
This function is in $L_2$ and we denote $\|\psi\|^2=\int_\R \psi^2(x) dx < 2$. 
Note that, since the function $\psi$ is infinitely differentiable with compact support,  we have:
\[
\lambda^* = \max_{0 < s <\smax} \sup_{x\neq y} \frac{\left|\psi^{(\lfloor s\rfloor)}(x) - \psi^{(\lfloor s\rfloor)}(y)\right|}{|x-y|^s} < +\infty.
\]

Now we consider $0<\alpha<1/8$,
$$\kappa_*^2=\frac{1}{32}\left(\frac{\|\psi\|^{2}}{2}\right)^{\! L}
c_1^2$$ and
$$c_1=\min\left\{ 1, M, (2L\lambda^*)^{-1}, \left(\frac{\alpha\log(2)}{8\cdot 2^L}\right)^{1/2} \right\}.$$


We consider the bandwidth
\[
    h = \left(\frac{1}{\Lambda_\ell^2\, n}\right)^{\frac1{2 s_\ell+\ell}}
\]
and we set $R=1/(2h)$. We assume, without loss of generality, that $R$ is an integer and $nh^\ell\ge 1$. Let $\mathcal{R}=\{0,\dotsc,R-1\}^\ell$ and define, for any $r=(r_1,\dotsc,r_\ell)\in\mathcal{R}$, the function $\phi_r:\Delta_\ell\to\R$ by:
\[
	\phi_r(y) = \prod_{i=1}^\ell \psi\left(\frac{y_i-x_i^{(r)}}{h}\right).
\]
where $x_i^{(r)}=(2r_i+1)h$. Finally, for any $w:\mathcal{R}\to\{0,1\}$ we define:
\[
	g_w = \rho_n\sum_{r\in\mathcal{R}} w(r) \phi_r
\]
where
\[
\rho_n = c_1 \Lambda_\ell h^{s_\ell} = c_1 \frac 1{\sqrt{nh^\ell}}.
\]

\paragraph{Proof of (ii).}
Set $w:\mathcal{R}\to\{0,1\}$. The following property can be readily verified:
\[
    \|g_w\|_{\Delta_\ell}^2 = |w|\,\|\psi\|^{2\ell}\rho_n^2 h^\ell
    \qquad\text{where}\qquad
    |w| = \left(\sum_{r\in\cR} w(r)\right) \leq R^\ell = \frac{1}{2^\ell h^\ell}.
\]
This implies that
\begin{equation}\label{eq:norme2gw}
    \|g_w\|_{\Delta_\ell}^2 \leq \left(\frac{\|\psi\|^{2}}{2}\right)^{\! \ell}
\rho_n^2\le c_1^2\frac{1}{nh^\ell}\le \ell! M^2.
\end{equation}
Moreover note that, for any $y\in\Delta_\ell$ and $\alpha=(\alpha_1,\dotsc,\alpha_\ell)$ such that $|\alpha|=\lfloor s_\ell\rfloor$, we have:
\begin{align*}
    D^\alpha\phi_r(y)
    &=
    \frac1{h^{|\alpha|}} \prod_{i=1}^\ell
    \psi^{(\alpha_i)}\left(\frac{y_i - x_i^{(r)}}{h}\right)
\end{align*}
which implies that, for any $z\in\Delta_\ell$ we have
\begin{align*}
\left|D^\alpha\phi_r(y)-D^\alpha\phi_r(z)\right|
&\leq  \frac{\|\psi\|_\infty^{\ell-1}}{h^{|\alpha|}}\sum_{i=1}^\ell
\left|
    \psi^{(\alpha_i)}\left(\frac{y_i-x_i^{(r)}}{h}\right) - \psi^{(\alpha_i)}\left(\frac{z_i-x_i^{(r)}}{h}\right)
\right| \\
&\leq \frac{\lambda^*}{h^{s_\ell}}\sum_{i=1}^\ell |y_i-z_i|^{s_\ell - \lfloor s_\ell\rfloor }\\
&\leq \frac{\ell\lambda^*}{h^{s_\ell}} |y-z|^{s_\ell - \lfloor s_\ell\rfloor}
\end{align*}
This also implies, since the function $\psi$ vanishes outside $(-1,1)$, that
\begin{align}
\left|D^\alpha g_w(y)-D^\alpha g_w(z)\right|
&\leq \left(2\ell\lambda^*\right)
\frac{\rho_n}{h^{s_\ell}}|y-z|^{s_\ell - \lfloor s_\ell\rfloor} \\
&\leq c_1\left(2\ell\lambda^*\right) \Lambda_\ell|y-z|^{s_\ell - \lfloor s_\ell\rfloor}\\
&\leq \Lambda_\ell|y-z|^{s_\ell - \lfloor s_\ell\rfloor}.\label{eq:diffgw}
\end{align}
Using \eqref{eq:diffgw}, we deduce that  $g_w$ belongs to $\cH_\ell(s_\ell, \Lambda_\ell)$. Combining with \eqref{eq:norme2gw}, (ii) is fulfilled.

\paragraph{Proof of (i) and (iii).} Using Lemma 2.9 of \citet{MR2724359}, there exists a set $\cW\subset\{w:\cR\to \{0,1\}\}$ such that the null function belongs to $\cW$, $\log_2 |\cW|\ge R^\ell/8$ and
\[
    \forall w\neq w'\in\cW , \ \sum_{r\in\mathcal{R}}\left| w(r) -w'(r)\right|\ge R^\ell/8.
\]
Let $w, w'\in \cW$ such that $w\neq w'$. We have
\begin{align}
    \|g_w-g_{w'}\|^2_{\Delta_\ell}=&\rho_n^2\sum_{r\in\cR}(w(r)-w'(r))^2 \|\phi_r\|_{\Delta_\ell}^2\\
    =&\rho_n^2\sum_{r\in\cR}|w(r)-w'(r)| h^{\ell}\|\psi\|^{2\ell}\\
    \ge & \rho_n^2h^{\ell}\|\psi\|^{2\ell} R^\ell/8\\
    \ge & \frac{1}{8}\left(\frac{\|\psi\|^{2}}{2}\right)^{\! L}
c_1^2\Lambda_\ell^2h^{2s_\ell}\\
\ge & 4 \kappa_*^2 \phi_n^2(s,\Lambda).
\end{align}
Then Assumptions (i) and (iii) are fulfilled.

\paragraph{Proof of (iv).}
Using \eqref{eq:norme2gw}, we deduce that using the definition of $c_1$
\begin{align}
    \frac{n}{|\mathcal{W}|}\sum_{\omega\in\mathcal{W}}\| g_{\omega}\|_{\Delta_\ell}^2\le & c_1^2h^{-\ell}\\
    \le & \frac{\alpha}{8}\left(2h\right)^{-\ell}\log(2)\\
    \le & \alpha \log |\cW|.
\end{align}
Then Assumption (iv) is fulfilled.

\subsection{Proof of Theorem~\ref{theo:2} }

We have using \eqref{eq:risk-decomposition} and~\eqref{eq:kappa0} that
\begin{align}
    \left( \e |\hat{m}_L(W)-m(W)|^p \right)^{1/p} 
    & \le \kappa_0n^{-1/2}+\sum_{\ell=1}^L \frac{\kc_\ell(p)}{\sqrt{\ell!}}
   \left(\e \left\|\hat f_{\hat{h}_\ell}^{(\ell)}-f_{\ell}\right\|^2_{\Delta_\ell}\right)^{1/2}\label{eq:controlmhat}.
\end{align}
 Let $\ell\in\{1,\ldots,L\}$. Let $h\in\fH_\ell$.
 We have
 \begin{align*}
     \left\|\hat f_{\hat{h}_\ell}^{(\ell)}-f_{\ell}\right\|_{\Delta_\ell}&\le  \left\|\hat f_{h}^{(\ell)}-f_{\ell}\right\|_{\Delta_\ell}+ \left\|\hat f_{\hat{h}_\ell\vee h}^{(\ell)}-\hat f_{ h}^{(\ell)}\right\|_{\Delta_\ell}+\left\|\hat f_{\hat{h}_\ell}^{(\ell)}-\hat f_{\hat{h}_\ell\vee h}^{(\ell)}\right\|_{\Delta_\ell}\\
      &\le \left\|\hat f_{h}^{(\ell)}-f_{\ell}\right\|_{\Delta_\ell}+ B(\ell,h)+M(\ell,\hat{h}_\ell)\\ &\qquad + B(\ell,\hat{h}_\ell)+ M(\ell,h)+ 2M(\ell,\hat{h}_\ell\vee h)\\
     &\le   \left\|\hat f_{h}^{(\ell)}-f_{\ell}\right\|_{\Delta_\ell} +4(B(\ell,h)+M(\ell,h)).
 \end{align*}
Then we have
\begin{align}
    \left(\e\left\|\hat f_{\hat{h}_\ell}^{(\ell)}-f_{\ell}\right\|^2_{\Delta_\ell}\right)^{1/2}&\le
    \left(\e\left\|\hat f_{h}^{(\ell)}-f_{\ell}\right\|^2_{\Delta_\ell}\right)^{1/2}+ 4 M(\ell,h) + 4\left( \e B^2(\ell,h) \right)^{1/2}.\label{eq:controlGL}
\end{align}
Note that we have
\begin{align}
  \left(\e\left\|\hat f_{h}^{(\ell)}-f_{\ell}\right\|^2_{\Delta_\ell}\right)^{1/2}
  &\le \left\|\e\hat f_{h}^{(\ell)}-f_{\ell}\right\|_{\Delta_\ell} + V(\ell)\nonumber\\
  &\le \left\|\e\hat f_{h}^{(\ell)}-f_{\ell}\right\|_{\Delta_\ell} + \frac{\sqrt{V_1(\ell)}+\sqrt{V_2(\ell)}}{\sqrt{n}}\nonumber\\
  &\le C\left( \left\|\e\hat f_{h}^{(\ell)}-f_{\ell}\right\|_{\Delta_\ell}+M(\ell,h)\right)\label{eq:controlris},
\end{align}
where we use the properties of $V_1(\ell)$ and $V_2(\ell)$ stated in page~\pageref{eq:V2}. In the following, we will demonstrate that
\begin{equation}\label{eq:controlBell}
    \left( \e B^2(\ell,h) \right)^{1/2}\le C\left(\max_{h'\le h} \left\|\e\hat f_{h'}^{(\ell)}-f_{\ell}\right\|_{\Delta_\ell}+ M(\ell,h)\right)+ O\left(\sqrt{\frac{\log n}{n}}\right).
\end{equation}

Combining \eqref{eq:controlGL} with \eqref{eq:controlris} and \eqref{eq:controlBell}, we obtain that
\[
    \left(\e\left\|\hat f_{\hat{h}_\ell}^{(\ell)}-f_{\ell}\right\|^2_{\Delta_\ell}\right)^{1/2}\le C\left(\max_{h'\le h} \left\|\e\hat f_{h'}^{(\ell)}-f_{\ell}\right\|_{\Delta_\ell}+ M(\ell,h)\right)+ O\left(\sqrt{\frac{\log n}{n}}\right).
\]
Theorem~\ref{theo:2} is then a direct consequence of the above inequality and~\eqref{eq:controlmhat}.

\paragraph{Proof of \eqref{eq:controlBell}}

Now let us control $B(\ell,h)$ for $h\in\fH_\ell$. We have
\begin{align*}
B(\ell,h)&=\max_{h'\le h}\left\{ \|\hat{f}_{h'}^{(\ell)}-\hat{f}_{h\vee h'}^{(\ell)}\|_{\Delta_\ell} -M(\ell,h') -M(\ell, h\vee h')  \right\}_+\\
& \le \max_{h'\le h}\left\{ \{\|\hat{f}_{h'}^{(\ell)}-\e\hat{f}_{ h'}^{(\ell)}\|_{\Delta_\ell} -M(\ell,h')\}_+\right\} + \{\|\hat{f}_{h}^{(\ell)}-\e\hat{f}_{ h}^{(\ell)}\|_{\Delta_\ell} -M(\ell,h)\}_+ \\
&\qquad +2\max_{h'\le h} \|\e\hat{f}_{ h'}^{(\ell)}-f_\ell\|_{\Delta_\ell}.
\end{align*}
Then
\begin{align}\label{eq:controlBell2}
    \left( \e B^2(\ell,h) \right)^{1/2}\le 2\max_{h'\le h} \|\e\hat{f}_{ h'}^{(\ell)}-f_\ell\|_{\Delta_\ell}+ 2\left(\#(\fH_\ell) \max_{h'\le h}A_\ell(h')   \right)^{1/2}
\end{align}
where
\[
    A_\ell (h') = \e\left[\left\{\left\|
    \hat f_{h'}^{(\ell)}-\e\hat f_{h'}^{(\ell)}
    \right\|_{\Delta_\ell}-M(\ell,h')\right\}_+^2\right].
\]
We have
\[
    \hat f_{h'}^{(\ell)}-\e\hat f_{h'}^{(\ell)}
      = \sum_{k=1}^L \fU_{k,\ell} + \fU_{\!\ell}
\]
where
\[
    k!\fU_{k,\ell}(t) =
     \frac1n \sum_{i=1}^n\bigg\{
    \Theta_{i,k} \xi_{i,\ell}(t,{h'})
    - \e\left(
    \Theta_{i,k} \xi_{i,\ell}(t,{h'})
    \right)
    \bigg\}
\]
and
\[
    \fU_{\!\ell}(t) =
    \frac1n\sum_{i=1}^n
    \varepsilon_i \xi_{i,\ell}(t,{h'}).
\]
Using these notations we have:
\[
    A_\ell({h'}) \leq 2^{L+1} \left(
    \sum_{k=0}^L \e\left\{\|\fU_{k,\ell}\|_{\Delta_\ell}-T(k,\ell)\right\}_+^2
    + \e\left\{\|\fU_{\ell}\|_{\Delta_\ell}-T(\ell)\right\}_+^2
    \right),
\]
where
\[
T(\ell)=\frac{\mu_4\left(1+4\sqrt{\log (1/(h')^\ell)}\right)}{\sqrt{n(h')^\ell}}
\]
and
\[
T(k, \ell)=\frac{\kc_4(k)Mk!\left(1+4\sqrt{\log (1/(h')^\ell)}\right)}{\sqrt{n(h')^\ell}}.
\]
Now note that
\[
    \left(\e\varepsilon_i^4\right)^{1/4}=\mu_4<\infty
\]
and  for $k\in\{0,\ldots,L\}$
\[
    \left(\e\Theta_{i,k}^4\right)^{1/4}\le \kc_4(k)Mk!.
\]
Using Lemma~\ref{lem:5} with $\fU=\fU_\ell$ (respectively $\fU=k! \fU_{k,\ell}$), $T=T(\ell)$  (respectively $T=T(k,\ell)$) and $\chi_i=\varepsilon_i$ (respectively $\chi=\Theta_{i,k}$), we deduce that for all for $h'\le h$
\[
    A_\ell(h')\le C n^{-1}.
\]
This implies that
\begin{equation}\label{eq:controlBell3}
\#(\fH_\ell) \max_{h'\le h}A_\ell(h') \le C (\log n)n^{-1}.
\end{equation}
Now \eqref{eq:controlBell2} and \eqref{eq:controlBell3} entail~\eqref{eq:controlBell}.

\subsection{Proof of Theorem~\ref{theo:3}.}

Let $s=(s_1,\ldots,s_L)\in [\smin,\smax)^L$, $\Lambda\in(0,+\infty)^L$, $M>0$ and $f\in \kM_{s,\Lambda, L, M}$. Define for $\ell\in\{1,\ldots,L\}$
\[
   k_\ell=\left\lfloor\frac{1}{2s_\ell+\ell}\log\left(\frac{n}{\log n}\right)\right\rfloor.
\]
For $n$ large enough, we have $h_\ell=e^{-k_\ell}\in\fH_\ell$. Using  \eqref{eq:bias-term} and \eqref{eq:defMell}, Theorem~\ref{theo:2} implies that
\begin{align*}
      R_2(\hat{m}_L, m)
      &\le C\sum_{\ell=1}^L \left\{\max_{h'\leq h_\ell} (h')^{s_\ell} + M(\ell, h_\ell)  \right\} + \Upsilon_2\left(\frac{\log n}{n }\right)^{1/2}\\
      &\leq  C\sum_{\ell=1}^L \left\{h_\ell^{s_\ell} + \left(\frac{\log n}{nh_\ell^\ell}\right)^{1/2}  \right\} + \Upsilon_2\left(\frac{\log n}{n}\right)^{1/2}\\
      &\leq C\tilde\phi_n(s) + \Upsilon_2\left(\frac{\log n}{n}\right)^{1/2}\\
      &\leq C\tilde\phi_n(s).
\end{align*}
where $C$ is a constant that changes from line to line and depends on $s$, $\Lambda$ and $M$. Since $C$ does not depend on $m\in\kM(s,\Lambda,L,M)$, this ends the proof.

\subsection{Proof of Lemma~\ref{lem:1}.}

We have
\begin{align}
    \left(\e |\xi_{\ell}(t,h)|^{2r}\right)^{1/r}&\le \kc_\ell^2(2r)\e |\xi_{\ell}(t,h)|^{2}\\
    &\le \kc_\ell^2(2r) \ell !\int_{\Delta_\ell} \left(K_{h}^{(\ell)}(t,u)\right)^2 \D u\\
    &\le \kc_\ell^2(2r) \ell ! 2^\ell \|\kkk\|_{[0,1]}^{2\ell}h^{-\ell}\\
    &\le b_{\ell,r}h^{-\ell}.
\end{align}

Moreover we have
\begin{align}
   \e\left(|\xi_{\ell}(t,h)|^r\1_{|\xi_{\ell}(t,h)| > \varphi}\right)
    &\leq \left(\e |\xi_{\ell}(t,h)|^{2r}\right)^{1/2} \left(\p\left(|\xi_{\ell}(t,h)| > \varphi\right)\right)^{1/2}\\
    &\leq (b_{\ell,r}h^{-\ell})^{r/2}\left(\frac{\e |\xi_{\ell}(t,h)|^{2q}}{\varphi^{2q}}\right)^{1/2}\\
    & \leq (b_{\ell,r})^{r/2} (b_{\ell,q})^{q/2} h^{-\ell(r+q)/2}\varphi^{-q}.
\end{align}

\subsection{Proof of Lemma~\ref{lem:5}}
{In this proof, $C$ is a positive constant that changes of value from line to line. }
Since $h$ is fixed, we simplify the notation and use in the proof $\xi_{i,\ell}(t,h)=\xi_{i,\ell}(t)$ and $\xi_{\ell}(t,h)=\xi_{\ell}(t)$. Now, we have for $k\ge 1$:
\begin{align}
\fU (t)    &= \bar\eta_0(t) + \eta_1(t) + \eta_2(t) + \eta_3(t)
\end{align}
where
\begin{align}
    \bar\eta_0(t) &= \eta_0(t) - \e\eta_0(t)
    \quad\text\quad
    \eta_0(t) =  \frac1n\sum_{i=1}^n \chi_i\1_{|\chi_i| \leq \varphi(n)} \xi_{i,\ell}(t) \1_{|\xi_{i,\ell}| \leq \psi(n)}\\
    \eta_1(t) &= \e\eta_0(t) - \e(\chi\xi_{\ell}(t))\\
    \eta_2(t) &= \frac1n\sum_{i=1}^n \chi_i\1_{|\chi_i| > \varphi(n)} \xi_{i,\ell}(t)\\
    \eta_3(t) &= \frac1n\sum_{i=1}^n \chi_i\1_{|\chi_i| \leq \varphi(n)} \xi_{i,\ell}(t) \1_{|\xi_{i,\ell}(t)| > \psi(n)}.
\end{align}
where for any $\ell=1,\ldots, L$:
\[
\psi(n)=n^\alpha
\quad\text{with}\quad
{\alpha = \alpha(\smin,\ell) =\frac{\smin+\ell}{4\smin+2\ell}}
\]
and
\[
\varphi(n) = n^\beta
\quad\text{with}\quad\beta = 1/2 - \alpha > 0.
\]
Note that both $\alpha$ and $\beta$ are positive numbers.
\subsubsection{Control of $\eta_1$}\label{sec:eta1}

We have
\begin{align}
\left|\eta_1(t)\right| &= \left|-\e \left( \chi\xi_\ell(t)\1_{|\chi|>\varphi(n)}\right)- \e \left( \chi\xi_\ell(t)\1_{|\chi|\leq \varphi(n)}\1_{|\xi_\ell(t)|> \psi(n)}\right)\right|\\
&\le \ \e \left( |\chi\xi_\ell(t)|\1_{|\chi|>\varphi(n)}\right)+ \e \left( |\chi\xi_\ell(t)|\1_{|\chi|\leq \varphi(n)}\1_{|\xi_\ell(t)|> \psi(n)}\right)\label{eq:controleta1}
\end{align}
Note that we have using Cauchy-Schwarz and Markov inequality
\begin{align}
   \e \left( |\chi\xi_\ell(t)|\1_{|\chi|>\varphi(n)}\right) &\le  \left(\e\left(|\chi|^2\1_{|\chi|>\varphi(n)}\right) \e |\xi_\ell(t)|^2\right)^{1/2}\\& \le {(b_{\ell,1})^{1/2} h^{-\ell/2}} a_2 \left(\p \left(|\chi|>\varphi(n)\right)  \right)^{1/4}\\
    &\le {(b_{\ell,1})^{1/2} h^{-\ell/2}}  a_2  \left(\e \left(|\chi|^{8/\beta}\right)\right)^{1/4}  \\
    &\le {(b_{\ell,1})^{1/2} h^{-\ell/2}}  a_2  (a_{4/\beta})^{2/\beta}  (\varphi(n))^{-2/\beta}\\
    &\le {Cn^{-1}}.\label{eq:controleta1bis}
\end{align}

Moreover since $h\in \fH_\ell$, using Lemma~\ref{lem:1} with $q=\frac{8\smin+6\ell}{\smin}$, we have
\begin{align}
    \left(\e \left( |\chi\xi_\ell(t)|\1_{|\chi|\leq \varphi(n)}\1_{|\xi_\ell(t)|> \psi(n)}\right)\right)^2\le & a_1^2 \e \left(|\xi_\ell(t)|^2\1_{|\xi_\ell(t)|> \psi(n)}\right)\\ &\le {a_1^2} (b_{\ell,2}) (b_{\ell,q})(\psi(n))^{-q}h^{-\ell(2+q)/2}\\
     &\le Cn^{-2}\label{eq:controleta1ter}.
\end{align}
Now using \eqref{eq:controleta1}, \eqref{eq:controleta1bis} and \eqref{eq:controleta1ter}, we finally obtain
\[
    \e\left(\|\eta_1\|_{\Delta_\ell}^2\right) \le Cn^{-1}.
\]

\subsubsection{Control of $\eta_2$}\label{sec:eta2}

Define
\[
    \mu_n(t) = \e\left(\chi\1_{|\chi| > \varphi(n)} \xi_\ell(t)\right).
\]
We have
\begin{align}
    \e\left(\|\eta_2\|_{\Delta_\ell}^2\right)
    &\leq 2(A + B)
\end{align}
where using Lemma~\ref{lem:1} with $q=4/\beta$
\begin{align}
    A & = \frac1n \int_{\Delta_\ell} \e\left(\chi \1_{|\chi| > \varphi(n)}\xi_\ell(t)-\mu_n(t)\right)^2 dt\\
    &\leq \frac1n
    \int_{\Delta_\ell}
    \left(
        \e\left(\chi^4 \1_{|\chi| > \varphi(n)}\right)
        \e\xi_\ell^4(t)
    \right)^{1/2} dt \\
    &\leq n^{-1} a_4^2 \p\left( |\chi|>\varphi(n) \right)^{1/4} b_{\ell,2}h^{-\ell}\\
    &\leq a_4^2 b_{\ell,2} \p\left( |\chi|>\varphi(n) \right)^{1/4}\\
    &\le a_4^2 b_{\ell,2} (a_{q/2})^{q/4}\left(\varphi(n)\right)^{-q/4}\\
    &\le Cn^{-1}
\end{align}
and following \eqref{eq:controleta1bis}
\[
    B
    = \frac1n \int_{\Delta_\ell} \mu_n^2(t) dt\le Cn^{-1}.
\]

\subsubsection{Control of $\eta_3$}\label{sec:eta3}

Define
\[
    \nu_n(t) = \e\left(\chi\1_{|\chi| \leq \varphi(n)} \xi_\ell(t)\1_{|\xi_\ell(t)| > \psi(n)}\right).
\]
We have
\begin{align}
    \e\left(\|\eta_3\|_{\Delta_\ell}^2\right)
    &\leq 2(A + B)
\end{align}
where
\[
    A = \frac1n \int_{\Delta_\ell} \e\left(\chi\1_{|\chi| \leq \varphi(n)} \xi_\ell(t)\1_{|\xi_\ell(t)| > \psi(n)}-\nu_n(t)\right)^2 dt
\]
and
\begin{align}
    B
    &= \frac1n \int_{\Delta_\ell} \nu_n^2(t) dt
\end{align}
Note that using similar arguments as above with $q=(8\smin+4\ell)/\smin$,
\begin{align}
    A
    &\leq \frac1n
    \int_{\Delta_\ell}
    \left(
        \e\left(\chi^4 \right)
        \e\left(\xi_\ell^4(t) \1_{|\xi_\ell(t)| > \psi(n)}\right)
    \right)^{1/2} dt \\
    &\leq a_2^2b_{\ell,4} (b_{\ell,q})^{q/4} (\psi(n))^{-q/2} h^{-\ell q/4}\\
    &\leq C n^{-1}
\end{align}
and following \eqref{eq:controleta1ter} $B\le Cn^{-1}$.

\subsubsection{Control of $\bar\eta_0$}\label{sec:eta0}

We have to bound
\begin{align}
    \e\left\{\|\bar\eta_{0}\|_{\Delta_\ell} - T \right\}_+^2
    &\leq \int_0^{+\infty} \p\left(\|\bar\eta_{0}\|_{\Delta_\ell} - T > \sqrt{u}\right) du\\
    &\leq 2\int_0^{+\infty} u\p\left(\|\bar\eta_{0}\|_{\Delta_\ell} > u + T \right) du.
\end{align}
Note that, using duality arguments, there exists  a countable set $\cS$ of functions $s\in\Leb^2(\Delta_\ell)$ such that $\|s\|_{\Delta_\ell}\leq 1$ and
\begin{align}
    \|\bar\eta_{0}\|_{\Delta_\ell}
    &= \sup_{s \in \cS} \int_{\Delta_\ell} s(t)\bar\eta_{0}(t) dt\\
    &= \frac{2\varphi(n)\psi(n)}n Z=\frac{2Z}{\sqrt{n}}
\end{align}
where
\[
    Z = \sup_{s\in\cS} \sum_{i=1}^n X_{i,s}
\]
and, for $s\in\cS$, we have:
\[
    X_{i,s} = \int_{\Delta_\ell} s(t)X_i(t) dt
    \qquad\text{with}\qquad
    X_i(t) = g_i(t)-\e g_i(t)
\]
and
\[
    g_i(t) = \frac{1}{2\sqrt{n}}
    \chi_i \1_{\{|\chi_i|\leq\varphi(n)\}}
    \xi_{i,\ell}(t)\1_{\{|\xi_{i,\ell}(t)|\leq\psi(n)\}}.
\]

Note that we have both $\e(X_{i,s}) = 0$ and $\|X_{i,s}\|_\infty \leq 1$. Now, let us control:
\[
    v = 2\e\left(\sup_{s\in\cS} \sum_{i=1}^n X_{i,s}\right) +  n\sup_{s\in\cS} \e X_{1,s}^2.
\]
Using Cauchy-schwarz's inequality and Fubini's theorem we obtain:
\begin{align}
    \e X_{1,s}^2
    &\leq \left(\int_{\Delta_\ell} s^2(t) dt \right) \left(\int_{\Delta_\ell} \e X_i^2(t) dt \right)\\
    &\leq \int_{\Delta_\ell} \e g_i^2(t) dt \\
\end{align}
We have
\begin{align}
    \int_{\Delta_\ell} \e g_i^2(t) dt
    &\leq \frac{a_2^2b_{\ell,2}}{4nh^{\ell}}
\end{align}
and
\begin{align}
    \e\left(\sup_{s\in\cS} \sum_{i=1}^n X_{i,s}\right)
    &= \e\left\|\sum_{i=1}^n X_i(\cdot)\right\|_{\Delta_\ell}\\
    &\leq \left(\e\left\|\sum_{i=1}^n X_i(\cdot)\right\|_{\Delta_\ell}^2\right)^{1/2}\\
    &\leq \frac{1}{4}\left(\int_{\Delta_\ell}\e \chi^2\xi^2_\ell(t)\D t\right)^{1/2}\\
    &\leq \frac{a_2\sqrt{b_{\ell,2}}}{4h^{\ell/2}}
\end{align}

Combining the previous results we have:
\[
    v \leq \theta_n+\sqrt{\theta_n}\quad  \text{with}\quad
    \theta_n = \frac{a_2^2b_{\ell,2}}{4h^{\ell}}.
\]
Define
\[
T = (1+\delta) \frac{a_2\sqrt{b_{\ell,2}}}{2\sqrt{nh^\ell}}=(1+\delta)\sqrt{\frac{\theta_n}{n}}.
\]
We have:
\begin{align}
    \p\left(\|\bar\eta_{0}\|_{\Delta_\ell} > u + T\right)
    &\le  \p\left(Z - \e Z > \frac{\sqrt{n}u}{2} + \frac{\delta}{2}\sqrt{\theta_n}\right).
\end{align}
Define:
\[
  \ka = \frac{\delta\sqrt{n\theta_n}}{2},
  \qquad
  \kb = \delta^2 \theta_n/4,
  \qquad
  \kc = \frac{\sqrt{n}}{3},
  \quad
  \kd = 2\theta_n + 2\sqrt{\theta_n} (1 +\delta/6).
\]
Using Bousquet's inequality we have:
\begin{align}
    \p\left(Z - \e Z > \frac{\sqrt{n}u}{2} + \frac{\delta}{2}\sqrt{\theta_n}\right)
    &\leq C_n(u) D_n(u)
\end{align}
where
\[
    C_{n}(u) = \exp\left(-\frac{nu^2}{4(\kc u+ \kd)}\right)
\]
and
\begin{equation}\label{eq:explicationlogn}
  D_n(u) = \exp\left(-\frac{\ka u + \kb}{\kc u + \kd}\right).
\end{equation}

Since $\ka\kd-\kb\kc>0$, we have, $D_n(u)\leq D_n(0)$, that is:
\[
    D_n(u) \leq \exp\left(-\frac{\kb}{\kd}\right)\leq \exp\left(-\frac{\delta^2\theta_n}{4(2\theta_n + 2\sqrt{\theta_n} (1 +\delta/6))}\right).
\]
Since $h\in \fH_\ell$, for $n$ large enough $(1+\delta/6) \leq \sqrt{\theta_n}$ we have $D_n(u) \leq h^\ell$. Moreover we have doing the change of variables $v=\sqrt{nh^\ell}u$
\begin{align}
  \e\left\{\|\bar\eta_{0}\|_{\Delta_\ell} - T\right\}_+^2
    &\leq 2\int_0^{+\infty} u\p\left(\|\bar\eta_{0}\|_{\Delta_\ell} > u + T\right) du\\
    &\leq 2h^\ell \int_0^{+\infty} u C_n(u) du\\
    &\leq \frac{C}{n} \int_0^{+\infty} v \exp\left(- \frac{Cv^2}{1+v}\right)\D v.
\end{align}
This implies that:
\begin{align}
  \e\left\{\|\bar\eta_{0}\|_{\Delta_\ell} - T\right\}_+^2
    &\leq C n^{-1}.\label{eq:eta0}
\end{align}

Combining results of Sections~\ref{sec:eta1}, \ref{sec:eta2}, \ref{sec:eta3} and \ref{sec:eta0}, we obtain \eqref{eq:conclusionlemma5}

\section*{Acknowledgements}
The authors have been supported by Fondecyt projects 1171335 and 1190801, and Mathamsud projects 19-MATH-06 and 20-MATH-05.

\bibliographystyle{plainnat}
\renewcommand*{\bibfont}{\small}
\bibliography{bibliography}

\end{document}